\documentclass{amsart}
\usepackage{import}
\usepackage{./preamble/Pakete}
\usepackage{./preamble/Befehle}
\bibliography{bibfinal.bib}

\begin{document}

\title[Universal Koszul duality for Kac--Moody groups]{Universal Koszul duality for Kac--Moody groups}

\author{Jens Niklas Eberhardt}
\email{mail@jenseberhardt.com}
\address{Johannes Gutenberg-Universität Mainz, Institut für Mathematik, Staudingerweg 9, 55128 Mainz, Germany}

\author{Arnaud Eteve}
\email{eteve@mpim-bonn.mpg.de}
\address{Max Planck Institute for Mathematics, Vivatsgasse 7, 53111 Bonn, Germany}

\begin{abstract}
We prove a monoidal equivalence, called \emph{universal Koszul duality}, between genuine equivariant $K$-motives on a Kac--Moody flag variety and constructible monodromic sheaves on its Langlands dual. The equivalence is obtained by a Soergel-theoretic description of both sides which extends results for finite-dimensional flag varieties by Taylor and the first author.

Universal Koszul duality bundles together a whole family of equivalences for each point of a maximal torus.
At the identity, it recovers an ungraded version of Beilinson--Ginzburg--Soergel and Bezrukavnikov--Yun's Koszul duality for equivariant and unipotently monodromic sheaves. 
It also generalizes Soergel-theoretic descriptions for 
monodromic categories on finite-dimensional flag varieties by Lusztig--Yun, Gouttard and the second author. For affine Kac--Moody groups, our work sheds new light on the conjectured quantum Satake equivalences by Cautis--Kamnitzer and Gaitsgory.

On our way, we establish foundations on six functors for reduced $K$-motives and introduce a formalism of constructible monodromic sheaves.
\end{abstract}
\maketitle

\section{Introduction}
\subsection{Motivation}
The seminal work of Beilinson--Ginzburg--Soergel \cite{beilinsonMixedCategoriesExtDuality,soergelKategorieMathcalPerverse1990, beilinsonKoszulDualityPatterns1996} and later refinements, for example \cite{soergelRelationIntersectionCohomology2000a,bezrukavnikovKoszulDualityKacMoody2013,soergelPerverseMotivesGraded2018,soergelEquivariantMotivesGeometric2018,
acharKoszulDualityKacMoody2019},
uncover a remarkable symmetry in the geometry and representation theory of a Kac--Moody group $G$ and its Langlands dual $\widehat{G}$.

The symmetry reveals itself as a monoidal equivalence, called Koszul duality, between \emph{equivariant} and \emph{unipotent monodromic} Hecke categories\footnote{For details, we refer to \cite{bezrukavnikovKoszulDualityKacMoody2013}, where the equivalence is shown for $\mathbb{Q}_\ell$-coefficients. Upcoming work by Ho and the first author based on \cite{hoRevisitingMixedGeometry2022,eberhardtMotivicSpringerTheory2022} will address the technicalities regarding semisimple Frobenius actions and prove an integral version of the equivalence in \emph{loc.~cit.}}
\begin{align*}
    \D_c^{\op{mix}}(B\bs G/B) \simeq \D^{\op{mix}}_{c}(\widehat{U}\bs \widehat{G}/ \widehat{U})_{\op{u-mon}}^\wedge
\end{align*}
which are certain mixed constructible sheaves on the stacks $B\bs G/B$ and $\widehat{U}\bs \widehat{G}/ \widehat{U}$, respectively. Here, $B\subset G$ (and $\widehat{U}\subset \widehat{G}$) denote the (unipotent radical of the) Borel subgroup. 

The natural actions on these two Hecke categories via Chern classes and monodromy, respectively, of the ring
$$R=\Z[\op{Pic}(\point/T)]=\Z[\pi_1(\widehat{T})]$$
factor through the completion $R_I^\wedge$ at the augmentation ideal $I$. So, Koszul duality can only see an \emph{infinitesimal} neighborhood of the point $1\in \Spec(R)=T.$

This raises the natural question of what is happening at other points $x\in T$. 
The goal of this article is to construct a \emph{universal Koszul duality} for Kac--Moody algebras, which naturally lives over $T$, and specializes to a family of dualities for any point $x\in T.$
This is, for example, motivated by a conjectural quantum $K$-theoretic Satake equivalence, see \Cref{sec:further}, and a conjectural universal Betti geometric Langlands, see \cite{ben-zviBettiGeometricLanglands2018}.
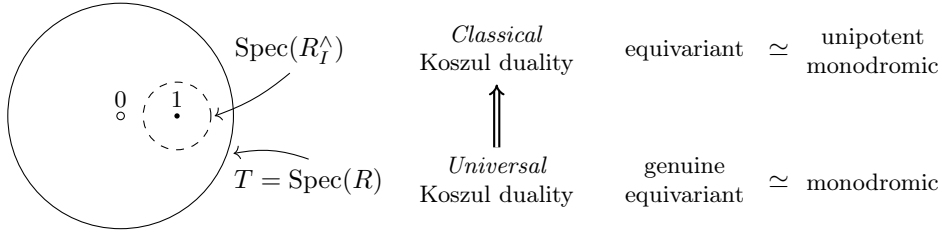
\begin{center}
    \begin{tikzpicture}[scale=0.5]
        \draw (0,0) circle (3cm);
        
        \draw (0,0) circle (0.1cm);
        \node[anchor=south, font=\small] at (0,0) {$0$};
        
        \coordinate (one) at (1.5,0);
        \fill (one) circle (2pt);
        \node[anchor=south, font=\small] at (one) {$1$};
        
        \draw[dashed] (one) circle (0.9cm);
        
        \node(spec_ri) at (4.5,1.75) {$\Spec(R_I^\wedge)$};
        \draw[->] (spec_ri.south) to[bend right=-10] (2.5,0);
        
        \node(spec_r) at (5,-1.75) {$T=\Spec(R)$};
        \draw[->] (spec_r.north) to[bend right=15] (3,-1);
        
        \node[text width=2.5cm, align=center, font=\small] (classical) at (10,1.75) {\emph{Classical}\\Koszul duality};

        \node[text width=2.5cm, align=center, font=\small] (eq) at (15,1.75) {equivariant};

        \node[text width=2.5cm, align=center, font=\small] (mon) at (20,1.75) {unipotent\\monodromic};

        \node[text width=2.5cm, align=center, font=\small] (universal) at (10,-1.75) {\emph{Universal}\\Koszul duality};
        
        \node[text width=2.5cm, align=center, font=\small] (geneq) at (15,-1.75) {genuine\\ equivariant};

        \node[text width=2.5cm, align=center, font=\small] (mono) at (20,-1.75) {monodromic};

        \draw[double equal sign distance, -Implies,thick] (universal) -- (classical);

        \draw[draw=none] (eq) to node[midway] {$\simeq$} (mon);

        \draw[draw=none] (geneq) to node[midway] {$\simeq$} (mono);
    \end{tikzpicture}
\end{center}

To achieve this, we follow a conjecture of the first author, see \cite{eberhardtMotivesKoszulDuality2022,eberhardtKtheorySoergelBimodules2024}, which proposes to replace equivariant sheaves with \emph{genuine equivariant $K$-motives} on the left and allow \emph{arbitrary monodromy} on the right. We will explain the two sides of the duality now.
\subsection{Genuine equivariant $K$-motives} 
The ring $R$ arises in a natural way as the \emph{genuine}\footnote{This is in contrast to equivariant formalisms defined via the Borel construction or \'etale descent where, by the Atiyah--Segal completion theorem, the action of $R$ factors through $R_I^\wedge$.} $T$-equivariant $K$-theory
$R=K_0^T(\point)=K_0(\point/T).$

Hence, we will replace equivariant constructible sheaves with a formalism that computes genuine equivariant $K$-theory: 
We consider the category of \emph{$K\!$-motives}
$$\DKbig(\Xx)=\Mod_{\KGL_\Xx}(\SH_{\op{mot}}(\Xx))$$
which is the category of modules of the $K$-theory spectrum in Hoyois' genuine equivariant stable motivic\footnote{Using motives is essential here. For $\op{KU}$-modules in the topological setting the shift $[2]$ acts as the identity and one could only obtain a $[2]$-periodic version of universal Koszul duality.} homotopy category, see \cite{hoyoisSixOperationsEquivariant2017}. $K$-motives are defined for `nice enough' stacks $\Xx$ and admit six functors, see \cite{hoyoisCdhDescentEquivariant2020, khanGeneralizedCohomologyTheories2022}.
For $\Xx$ smooth, the mapping space of the monoidal unit is the algebraic $K$-theory spectrum
$K(\Xx).$ 

For our purposes, the higher algebraic $K$-theory of the base field $k$ is irrelevant. We will hence work with \emph{reduced} $K$-motives $\DK(\Xx)$ which arise by modding out $K_{>0}(\Spec(k))$ from $\DKbig(\Xx)$, see \cite{eberhardtIntegralMotivicSheaves2023,eberhardtKmotivesSpringerTheory2024}. 
Most importantly, the ring $R$ arises as mapping space in the category of reduced $K$-motives on $\point/T$, 
$$R=\Map_{\DK(\point/T)}(\un,\un).$$ We discuss $K$-motives and an extension to so-called linearly reductive ind-pro-stacks in \Cref{sec:kmotives}.
\subsection{$K$-theoretic Hecke category} 
We assume that the Kac--Moody datum $\Dd$ of $G$ is free and of simply-connected type.
For example, we can take $G$ as the loop group of a simply-connected quasi-simple algebraic group extended by loop rotations. Basic results on Kac--Moody groups are recalled in \Cref{sec:kacmoody}.

In \Cref{sec:soergelktheoretichecke}, we define the \emph{$K$-theoretic Hecke category}
$$\HeckeK_{\Dd}=\DKbig_{\red}(B\bs G/B)_{\op{l.c.}}\subset \DK(B\bs G/B)$$
as a full subcategory of reduced $K$-motives on Hecke stack $B\bs G/B$ that are `locally constant' along the Bruhat cells $BwB.$ The main goal of \Cref{sec:soergelktheoretichecke} is to obtain a `Soergel-theoretic' description of $\HeckeK_{\Dd}$.

Multiplication on $G$ induces two monoidal structures $\conv$ and $\conv^!$ on $\HeckeK_{\Dd}$ defined via $*$- and $!$-functors, respectively. A first important result, see \Cref{thm:rigiditykmotivichecke}, is that the monoidal structures are equivalent and $\HeckeK_{\Dd}$ is rigid. Here, standard techniques do not apply since $G/B$ is not necessarily smooth and the six functors for $K$-motives only work for representable maps. We circumvent these issues by using results of Boyarchenko--Drinfeld \cite{boyarchenkoDualityFormalismSpirit2013}.

In analogy to Soergel's functor $\mathbb{H}$, we define a functor
$$\Kyp: \HeckeK_{\Dd}\to \D(R\otimes R)$$
which, in essence, maps a $K$-motive to its $K$-theory as a bimodule over $R.$ We show that the functor $\Kyp$ is monoidal, see \Cref{thm:monoidalityofkyp}, and prove an analog of Soergel's Erweiterungssatz: the functor $\Kyp$ is fully faithful when restricted to pure objects, see \Cref{thm:erweiterungssatz}. Moreover, we prove a formality result for pure objects, see \Cref{cor:formalityofkmotivichecke}, which relies on the fact that we removed the higher algebraic $K$-theory of the base field.

The essential image of pure objects under the functor $\Kyp$ is the category of $K$-theory Soergel bimodules $\SBimK_{\Dd}$ defined in \cite{eberhardtKtheorySoergelBimodules2024}. These are direct summands of the bimodules 
$$R\otimes_{R^{s_1}}\dots \otimes_{R^{s_n}} R$$
which arise from the equivariant $K$-theory of Bott--Samelson resolutions. 

Combining these results we arrive at the following Soergel-theoretic description.
\begin{theorem*}[\Cref{thm:soergelktheoreticheckecategory}]
    There is an equivalence of monoidal categories
    $$\HeckeK_{\Dd}\to \Ch^b(\SBimK_{\Dd})$$
    between the $K$-theoretic Hecke category and the category\footnote{Throughout this paper, we work with $\infty$-categories. On the level of homotopy categories, the target of the equivalence is the bounded homotopy category of chain complexes.} of bounded chain complexes of $K$-theoretic Soergel bimodules associated to $\mathcal{D}.$ 
\end{theorem*}

We note that for finite flag varieties, with rational coefficients and without the monoidal structure this result was shown in \cite{eberhardtKtheorySoergelBimodules2024}.

\subsection{Monodromic sheaves}
Dually to the $K$-theoretic side, the ring $R$ arises as the group ring of the fundamental group of the Langlands dual torus $\Z[\pi_1(\widehat{T})]$. The category $\D(R)$ is equivalent to the category of monodromic sheaves in $\D(\widehat{T}, \Z)$, that is, sheaves that are locally constant. Under this equivalence, $R$ corresponds to the free local system on $\widehat{T}.$ However, the free local system is infinite-dimensional and hence \emph{not constructible} as a $\Z$-sheaf. This is problematic since, a priori, this prohibits to use the important deep results from the theory of constructible sheaves. The goal of \Cref{sec:monodromicsheaves} is to address this issue and to develop a nice theory of \emph{constructible monodromic} sheaves.

The idea of monodromic sheaves is not new in the literature and exists in several sheaf-theoretic contexts. For étale sheaves, they were first considered by Verdier \cite{AST_1983__101-102__332_0} and defined, for a finite ring of coefficients $\Lambda$, to be the full subcategory of locally constant sheaves in $\D^b_c(\widehat{T}, \Lambda)$. This category, while easy to define, does not contain the free local system. To circumvent these issues, Bezrukavnikov-Yun \cite[Appendix A]{bezrukavnikovKoszulDualityKacMoody2013} introduce the category of unipotently monodromic sheaves on $\widehat{T}$ which is defined so that it satisfies
\begin{equation*}
\D_c(T)_{\op{u-mon}}^{\wedge} = \D^b(R^{\wedge}_I),
\end{equation*}
where $R^{\wedge}_I$ is the completion along the augmentation ideal. 

Our formalism of constructible monodromic sheaves on a complex algebraic variety $X$ with a $\widehat{T}$-action will use three equivalent definitions of the category of monodromic sheaves, each with their own advantages.

First, one may simply define $\D(X,\Z)_{\op{mon}}=\D(X,\Z)_{\op{l.c.}}$ as sheaves that are locally constant along the $T$-orbits. Second, one can use that the exponential map $\widehat{\mathfrak{t}} \to \widehat{T}$ is a contractible universal cover and define $\D(X,\Z)_{\op{mon}}=\D(X/\widehat{\mathfrak{t}} ,\Z)$ as $\widehat{\mathfrak{t}} $-equivariant sheaves. Third, we will introduce the following new definition of monodromic sheaves, inspired by \cite{gabberFaisceauxPerversEll1996}. We denote by $\Ll_{\widehat{T}}$ the rank one free local system on $\widehat{T}$. We observe that $\Ll_{\widehat{T}}$ is multiplicative (\Cref{lem:MultiplicativityFreeLocalSystem}), that is, there is an isomorphism
$$m^*\Ll_{\widehat{T}} = \Ll_{\widehat{T}} \boxtimes_R \Ll_{\widehat{T}}.$$
Once this property is established, we can define $\D(X)_{\op{mon}}=\D(X/(\widehat{T}, \Ll_{\widehat{T}}),R)$ as twisted $(\widehat{T}, \Ll_{\widehat{T}})$-equivariant sheaves with coefficients in $R$, see \cite{gaitsgoryLocalGlobalVersions2020}. 
We will show the following comparison result.
\begin{theorem*}[\Cref{lemCharacterizationMonodromic}]\label{thm:definitionMonoIntro} There is an equivalence between the categories $$\D(X,\Z)_{\op{l.c.}}, \D(X/\,\widehat{\mathfrak{t}} ,\Z)\text{ and }\D(X/(\widehat{T}, \Ll_{\widehat{T}}),R).$$ All three categories embed fully faithfully in the category $\D(X,\Z)$ with the same essential image we denote by $\D(X,\Z)_{\op{mon}}.$
\end{theorem*}
A similar approach was considered in the second author's thesis \cite{eteveMonodromiquesSheavesDeligneLusztig2023} in the étale context. Based on this, we propose the following definition.
\begin{definition*}[\Cref{def:monodromicSheaves}]
The category of \emph{constructible monodromic sheaves} $$\D_c(X,\Z)_{\op{mon}}\subset \D(X,\Z)_{\op{mon}}$$ 
is the full subcategory of monodromic sheaves which are constructible as $R$-sheaves, when seen as objects in $\D(X/(\widehat{T}, \Ll_{\widehat{T}}),R).$
\end{definition*}

This definition has several advantages.
It is still equipped with a $4$-functor formalism, see \Cref{lemmaFunctorialityMonodromic}. Moreover, it has a good notion of duality, see \Cref{subsectionDuality}.
Also, there is no restriction on the space $X$, contrary to \cite{bezrukavnikovKoszulDualityKacMoody2013} which requires the $\widehat{T}$-action to be free.

\subsection{Monodromic Hecke category} We assume that the root datum $\widehat{\Dd}$ of the Kac--Moody group $\widehat{G}$ is cofree and of adjoint type.\footnote{For ease of notation, the group $\widehat{G}$ will be denoted by $G$ in \Cref{sec:soergelmonodromic}.}

In \Cref{sec:soergelmonodromic}, we consider the \emph{universal monodromic Hecke category} $$\Heckemon_{\widehat{\Dd}}=\D_c(\widehat{U}\bs \widehat{G}/\widehat{U})_\mon$$
which is the category of constructible monodromic sheaves on the (enhanced) Hecke stack $\widehat{U}\bs \widehat{G}/\widehat{U}$. In particular, we do not impose any generalized character on the monodromy. 
The main goal of \Cref{sec:soergelmonodromic} is a Soergel-theoretic description of $\Heckemon_{\widehat{\Dd}}$. We also refer to \cite[Remark 4.11]{ben-zviBettiGeometricLanglands2018}, where a conjectural Bezrukavnikov-type description of monodromic sheaves on affine flag varieties is given.

The universal monodromic Hecke category has a perverse $t$-structure. There are perverse (co-)standard objects $\Delta_w,\nabla_w\in \Heckemon_{\widehat{\Dd}}$ which arise from the $!$- and $*$-pushforward of the free local system on each Bruhat cell. We consider tilting perverse sheaves which admit both a standard and costandard flag. Following \cite{beilinsonTiltingExercises2004a}, the Hecke category is equivalent to bounded chain complexes of tilting perverse sheaves.

There are two monoidal structures $\conv$ and $\conv^!$ on the Hecke category. 
Using rank one calculations of Taylor \cite{taylorUniversalMonodromicTilting2023}, we show that the (co-)standard and tilting objects behave as expected with respect to both monoidal structures. This and results by Boyarchenko--Drinfeld \cite{boyarchenkoDualityFormalismSpirit2013} allow us to show that both products are equivalent and that $\Heckemon_{\widehat{\Dd}}$ is rigid, see \Cref{thm:rigiditymonodromic}.

Following Li--Nadler--Yun \cite{liFunctionsCommutingStack2023}, we define the functor
$$\mathbb{V}:\Heckemon_{{\widehat{\Dd}}}\to \D(R\otimes R), A\mapsto 1^*\phi_\chi(i^*A)$$
as a stalk of the vanishing cycles along a character $\chi: \widehat{U}^-\to \mathbb{G}_a$ together with its left and right monodromy action. By \cite{taylorUniversalMonodromicTilting2023}, the functor is monoidal and sends the tilting sheaf $T_s$ supported on the minimal parabolic $P_s$ to the bimodule $R\otimes_{R_s}R.$
In \Cref{thmStrukturSatzMonodromic} we prove an analog of Soergel's Struktursatz, namely that $\mathbb{V}$ is fully faithful on tilting objects.

Combining these results we arrive at the following Soergel-theoretic description.
\begin{theorem*}[\Cref{thm:soergelmonodromicheckecategory}]
    There is an equivalence of monoidal categories
    $$\Heckemon_{\widehat{\Dd}}\to \Ch^b(\SBimK_{\Dd})$$
    between the universal monodromic Hecke category associated to $\widehat{\Dd}$ and the category of bounded chain complexes of $K$-theoretic Soergel bimodules associated to $\Dd.$ 
\end{theorem*}
For finite flag varieties this result was shown in \cite{taylorUniversalMonodromicTilting2023} using localization techniques.

\subsection{Universal Koszul duality} By combining the Soergel-theoretic descriptions of both Hecke categories, we obtain our main theorem.
\begin{theorem*}[Universal Koszul Duality, \Cref{thm:main}] There is a monoidal equivalence of categories
    $$\HeckeK_{\mathcal{D}}\simeq \Heckemon_{\widehat{\mathcal{D}}}$$
    between the $K$-theoretic and universal monodromic Hecke categories associated to Langlands dual Kac--Moody data. The equivalence exchanges pure $K$-motives and perverse monodromic tilting objects.
\end{theorem*}
\subsection{Specialization}
We now explain the relation of universal and classical Koszul duality.
Recall that $I$ denotes the augmentation ideal of $R.$ By completing at $I$ and tensoring with any coefficient ring $\Lambda$, universal Koszul duality specializes to a diagram of the form 
\[\begin{tikzcd}
	{\HeckeK_{\Dd}} & {\Heckemon_{\widehat{\Dd}}} \\
	{(\HeckeK_{\Dd})_{I}^\wedge\otimes \Lambda=\DKbig^{\acute{e}t}_{\red,c}(B\bs G/B,\Lambda)_{\op{l.c.}}} & {\D_c(\widehat{U}\bs \widehat{G}/\widehat{U},\Lambda)^\wedge_{\operatorname{u-mon}}=(\Heckemon_{\widehat{\Dd}})_I^\wedge\otimes \Lambda} \\
	{\D^{\op{mix}}_c(B\bs G/B,\Q)} & {\D^{\op{mix}}_c(\widehat{U}\bs \widehat{G}/\widehat{U},\Q)^\wedge_{\operatorname{u-mon}}}
	\arrow["\simeq", from=1-1, to=1-2]
	\arrow[from=1-1, to=2-1]
	\arrow[from=1-2, to=2-2]
	\arrow["\simeq", from=2-1, to=2-2]
	\arrow["\iota", dashed, from=3-1, to=2-1]
	\arrow["\simeq", from=3-1, to=3-2]
	\arrow["v"', dashed, from=3-2, to=2-2]
\end{tikzcd}\]
where by $\DKbig^{\acute{e}t}_{\red}$ we denote reduced $K$-motives with enforced \'etale descent.

The middle equivalence can be seen as an ungraded version of classical Koszul duality where the 
dashed arrows $v,\iota$ exist if $\Q\subset \Lambda$ and `forget' the extra grading induced by the mixed structure: $\iota(A(1)[2])=\iota(A)$ and $v(A(1))=v(A)$, see \cite{eberhardtMotivesKoszulDuality2022}. The bottom equivalence is classical Koszul duality.

If $\widehat{G}$ is reductive, \Cref{thm:soergelmonodromicheckecategory} implies the Soergel-theoretic description of the category of unipotently monodromic sheaves from Bezrukvanikov--Riche \cite{bezrukavnikovTopologicalApproachSoergel2020}. We can also specialize at any other ideal $\mathfrak{m}\subset R$ and recover the (completed) categories of monodromic sheaves on which the action of $\mathfrak{m}$ is locally finite
$$(\Heckemon_{\widehat{\Dd}})^\wedge_\mathfrak{m}\otimes \Lambda=\D_c(\widehat{U}\bs \widehat{G}/\widehat{U},\Lambda)^\wedge_{\mathfrak{m}\!\operatorname{-mon}}.$$
\Cref{thm:soergelmonodromicheckecategory} yields a Soergel-theoretic description for these categories. If $\widehat{G}$ is reductive, this recovers results of Lusztig--Yun \cite{lusztigEndoscopyHeckeCategories2020}, Gouttard \cite{gouttardPerverseMonodromicSheaves2021} and the second author \cite{eteveMonodromiquesSheavesDeligneLusztig2023}.
\subsection{Further directions}\label{sec:further} This paper is part of an ongoing program extending results in geometric representation theory/geometric Langlands to $K$-theory. We discuss some further directions.
\subsubsection{Parabolic/Whittaker duality} 
The universal Koszul duality as discussed should also admit a parabolic/Whittaker version, similar to the results in Bezrukavnikov--Yun \cite{bezrukavnikovKoszulDualityKacMoody2013}. 

On the equivariant side $B\bs G/B$ can be replaced by a quotient $P\bs G/Q$ for finitary parabolic subgroups $P,Q\subset G.$ 
By an argument as in \cite[Corollary III.6.10]{soergelEquivariantMotivesGeometric2018}, the resulting parabolic $K$-theoretic Hecke category admits a combinatorial description in terms of singular $K$-theory Soergel bimodules over $R^{W_P}\otimes R^{W_Q}$. 

It is reasonable to expect that the parabolic $K$-theoretic Hecke category for $P\bs G/Q$ should correspond to certain `bi-Whittaker' monodromic sheaves on $\widehat{U}\bs \widehat{G}/\widehat{U}$. However, handling (bi-)Whittaker objects for universal monodromic Betti sheaves seems to be subtle. A potential approach is to use the Kirillov model considered in \cite[Section 1.6]{gaitsgoryLocalGlobalVersions2020}. In the context of $D$-modules \cite{chenLanglandsDualRealization2023} gives a description of sheaves on $G$ with either equivariant, Whittaker or bi-Whittaker conditions, see also \cite{campbellAffineHarishChandraBimodules2021}.
\subsubsection{Quantum Satake equivalence} An interesting example is the parabolic $K$-theoretic Hecke category for affine Grassmannian of a reductive group $G$, equivariant with respect to the positive loop group and loop rotations
$$\DKbig_{\red}((G(k[[t]])\times \Gm)\bs\op{Gr}_{G})_{\op{l.c.}}.$$
By a conjecture of Cautis--Kamnitzer \cite{cautisQuantumKtheoreticGeometric2018}, there should be a quantum $K$-theoretic derived Satake equivalence 
$$\DKbig_{\red}((G(k[[t]])\times \Gm)\bs\op{Gr}_{G})_{\op{l.c.}}
\simeq 
\D^b_{U_q(\widehat{\mathfrak{g}})}(\mathcal{O}_q(\widehat{G}))$$
with a category of representations of the Langlands dual quantum group. By our work, the left-hand side admits a description in terms of maximally singular $K$-theory Soergel bimodules. It would be very desirable to achieve the same for the right-hand side.

Moreover, the universal parabolic/Whittaker Koszul duality provides a bridge between Cautis--Kamnitzer's and Gaitsgory's \cite{gaitsgoryTwistedWhittakerModel2007} approaches to a quantum Satake. We will discuss this elsewhere.

\subsection{Structure} In \Cref{sec:kacmoody} we record basic definitions and properties of Kac--Moody root data, groups and flag varieties, as well as the relation to loop groups of reductive groups. Moreover, we prove some results on $K$-theory Soergel bimodules for Kac--Moody groups.

\Cref{sec:kmotives} recalls the definition of (reduced) $K$-motives for linearly reductive stacks. We explain how to extend the formalism to ind-pro-stacks. In \Cref{sec:soergelktheoretichecke}, we define the $K$-theoretic Hecke category and prove a description in terms of $K$-theory Soergel bimodules.

In \Cref{sec:monodromicsheaves} we introduce constructible monodromic sheaves.  \Cref{sec:soergelmonodromic} defines the universal monodromic Hecke category and provides a description in terms of $K$-theory Soergel bimodules.

Finally, in \Cref{sec:maintheorem} we show the universal Koszul duality.

\subsection{Notations}
Denote by $\mathrm{Pr}$ the $\infty$-category of presentable stable $\infty$-categories with left adjoint functors, as defined in \cite[Section 5.5]{lurieHigherToposTheory2009}. For a ring $\Lambda$, denote by $\mathrm{Pr}_\Lambda$ the category of $\Lambda$-linear presentable stable categories. In a stable $\infty$-category $\mathcal{C}$, we denote the mapping spectrum by $\Map_\mathcal{C}$ and the set of homomorphisms by $\Hom_{\mathcal{C}}=\pi_0\Map_{\mathcal{C}}$.

Given a property (P) of morphisms in a category $\Cc$, a morphism $f$ of pro-objects in $\Cc$ is called pro-(P) if there is a presentation $f=\lim f_i$ such that each $f_i$ in $\Cc$ has property (P). We use the same convention for morphisms of ind- and ind-pro-objects.
\subsection{Acknowledgements} We thank Quoc Ho, Marc Hoyois, Adeel Khan and Peter Scholze for helpful discussions. We thank Wolfgang Soergel and Gurbir Dhillon for valuable comments on the first draft. The first author was supported by Deutsche Forschungsgemeinschaft (DFG), project number 45744154, Equivariant K-motives and Koszul duality. The second author was supported by the Max Planck Institute for Mathematics. 

\setcounter{tocdepth}{1} 
\tableofcontents
\section{Kac--Moody groups and $K$-theory Soergel bimodules}\label{sec:kacmoody}
We recall some basic properties of Kac--Moody root data, their associated groups and flag varieties as well as $K$-theory Soergel bimodules. For the construction of Kac--Moody groups, we follow \cite{mathieuConstructionGroupeKacMoody1989} and \cite{rousseauGroupesKacMoodyDeployes2016}. 
Fix a base field $k$. 
\subsection{Kac--Moody root data}\label{sec:kacmoodyrootdatum}
Let $\mathcal{D}=(X, \{\alpha_i\}_{i\in I}, \{\alpha_i^{\vee}\}_{i\in I})$ be a Kac--Moody root datum with generalized Cartan matrix $A = (a_{i,j})_{i,j \in I}.$ This means that $X$ is a finite-rank lattice, $\alpha_i\in X$, $\alpha_i^\vee\in X^\vee$ and $a_{i,j} = \langle \alpha_j, \alpha_i^{\vee}\rangle$ where by $\langle -,-\rangle$ we denote the evaluation pairing and $X^\vee$ the dual space of $X$. 

By dualizing and exchanging roots and coroots, one obtains the \emph{Langlands dual} Kac--Moody root datum $\LD{\mathcal{D}}=(X^\vee, \{\alpha_i^{\vee}\}_{i\in I}, \{\alpha_i\}_{i\in I})$
with Cartan matrix $A^{\operatorname{tr}}.$ 

The \emph{Weyl group} $W(\mathcal{D})\subset \op{GL}(X)$ associated to the Kac--Moody datum is the group generated by the \emph{simple reflections} $s_i$ with 
$$s_i(\lambda)=\lambda-\langle \lambda,\alpha_i^\vee \rangle \alpha_i \text{ for }\lambda\in X.$$

We call a Kac--Moody root datum $\mathcal{D}$ \emph{free} (or \emph{cofree}) if $\{\alpha_i\}_{i\in I}$ (or $\{\alpha^\vee_i\}_{i\in I}$) is linearly independent.
We say that $\mathcal{D}$ is of simply-connected type (or adjoint type) if for all $i\in I$, $\langle -,\alpha_i^\vee\rangle: X\to \Z$ (or $\langle \alpha_i,-\rangle: X^\vee\to \Z$) is surjective.

\subsection{Kac--Moody groups} \label{sec:kacmoodygroups}
To a Kac--Moody root datum $\mathcal{D}$ one can associate a Kac--Moody group, Borel subgroup and maximal torus
$G\supset B\supset T$ over $k$. Here, $G$ and $B$ are group ind-schemes and $T$ is the split torus with (co-)character lattice $X(T)=X$ and $Y(T)=X^\vee$, respectively.
We denote by $\mathfrak g\supset \mathfrak b \supset \mathfrak h$ the corresponding Lie algebras, by $e_i,f_i\in \mathfrak g$ the Chevalley generators and by $\Delta_+\subset\Delta$ the set of (positive) roots. 

We record some important properties of the Kac--Moody group.
The \emph{Weyl group} of the Kac--Moody group is $W(G,T)=N_G(T)/T$. The action of $W(G,T)$ on $X$ yields a surjection
$$W(G,T)\twoheadrightarrow W(\mathcal{D})\subset \op{GL}(X)$$
which is an isomorphism if the root datum is free. In this case, we simply write $W=W(G,T)= W(\mathcal{D})$.
For all $i\in I$, there is a map 
$\varphi_i:\op{SL}_2\to G$
such that $\varphi_i(\op{diag}(t,t^{-1}))=\alpha_i^\vee(t)\in T$ and 
$$d\varphi_i\left(\begin{pmatrix}
        0 & 1 \\
        0 & 0
        \end{pmatrix}\right)=e_i \text{ and }d\varphi_i\left(\begin{pmatrix}
            0 & 0 \\
            1 & 0
            \end{pmatrix}\right)=f_i.
$$
Denote the image of $\varphi_i$ by $G_i$. The map $\varphi_i$ is injective 
if and only if $\langle -,\alpha_i^\vee\rangle: X\to \Z$ is surjective. This motivates our definition of root data of simply-connected type in Section \ref{sec:kacmoodyrootdatum}.

The Lie algebra $\mathfrak n_+$, generated by the $e_i$, has a filtration by the ideals
$$\mathfrak n_+(m)=\prod_{\substack{\alpha \in \Delta_+\\ \op{ht}(\alpha)\geq m}} \mathfrak g_\alpha.$$
The completion with respect to the filtration is a pro-Lie algebra $\hat{\mathfrak n}_+$. Denote the associated pro-group by $U.$ Then $U$ has a filtration by normal subgroups $U(m)$ with $\op{Lie}(U(m))=\hat{\mathfrak n}_+(m)$ where $U/U(m)$ is a (finite dimensional) linear algebraic group and 
$U=\lim_n U/U(n).$
The Borel subgroup is a semidirect product $B=T\ltimes U$. 

For each subset $J\subset I$, there is a \emph{standard parabolic subgroup} $P_J$ generated by $B$ and $G_j$ for $j\in J.$ There is a Levi decomposition $P_J=L_J\ltimes U^J$ such that $L_J$ is the Kac--Moody group with Kac--Moody datum $\mathcal{D}_J=(X, \{\alpha_j\}_{j\in J}, \{\alpha_j^{\vee}\}_{j\in J})$ and $U^J\subset U.$

\subsection{Kac--Moody flag varieties} \label{sec:kacmoodyflagvarieties}
The \emph{Kac--Moody flag variety} is the quotient $G/B$. For $w\in W$ denote the \emph{Bruhat cell} $(G/B)_w=BwB/B\cong \mathbb{A}^{\ell(w)}$. Then $G/B$ is an ind-variety filtered by closed projective subvarieties
$$(G/B)_{\leq n}=\bigcup_{\substack{w\in W\\ \ell(w)\leq n}} (G/B)_{w}$$
for $n\geq 0.$ The \emph{Hecke stack} associated to the Kac--Moody datum $\mathcal{D}$ is the double quotient $$B\backslash G/B.$$ We now show that this is a \niceindprostack ind-pro stack. See \Cref{sec:extensiontoindstacks} for the definition of \niceindprostack (ind-)pro stacks.

Note that $wB/B$ is stabilised by $U\cap wUw^{-1}\supset U(m(w))$ where 
$$m(w)=\max\setbuild{\op{ht}(\alpha)}{\alpha\in \Delta_+\cap w^{-1}(\Delta_-)}+1.$$
Note that $m(w)$ is finite, since $|\Delta_+\cap w^{-1}(\Delta_-)|=\ell(w).$ For $n\geq 0$, let 
$$m(n)=\max\setbuild{m(w)}{w\in W, \ell(w)\leq n}$$
which is again finite.
Then $U(m(n))$ stabilises all points $wB/B$ for $\ell(w)\leq n.$ Since $U(m(n))\subset B$ is a normal subgroup, it acts trivially on $(G/B)_{\leq n}.$
Hence the action of $B$ on $(G/B)_{\leq n}$ factors through the finite dimensional group $B/U(m(n))$. Moreover, $\lim_{m\geq m(n)} B/U(m)= B$ and we see that 
$$(B\backslash G/B)_{\leq n}=\lim_{m\geq m(n)} \left(B/U(m)\backslash (G/B)_{\leq n}\right)$$
is a \niceprostack pro-stack whose transition maps are affine space bundles with fiber $U(m)/U(m+1)$.
The morphisms $(B\backslash G/B)_{\leq n}\to (B\backslash G/B)_{\leq n+1}$ are pro-closed immersions and we see that
$$B\backslash G/B=\colim_n (B\backslash G/B)_{\leq n}$$
has the structure of a \niceindprostack ind-pro stack.
\subsection{Affine Kac--Moody groups}
Let $\mathring{\mathcal{D}}=(\mathring{X}, \{\alpha_i\}_{i\in \mathring{I}}, \{\alpha_i^{\vee}\}_{i\in \mathring{I}})$ be a classical root datum associated to an indecomposable Cartan matrix $\mathring{A}$ for $\mathring{I}=\{1,\dots,\ell\}.$ Associated to this, we obtain a reductive group $\mathring{G}\supset \mathring{B}\supset \mathring{T}$ with Borel subgroup and maximal torus. 

The root datum $\mathring{\mathcal{D}}$ is free and cofree. Moreover, it is of simply-connected type if the derived subgroup $[\mathring G,\mathring G]$ is simply-connected and of adjoint type if the center $Z(\mathring G)$ is connected.

Let $\alpha_0=-\theta\in \mathring X$ and $\alpha_0^\vee=-\theta^\vee\in \mathring X^\vee$, 
where $\theta\in \mathring \Delta$ 
is the highest root in the root system of $\mathring{\mathcal D}.$ We obtain the Kac--Moody root datum
$\mathcal D=(\mathring X, \{\alpha_i\}_{i\in I}, \{\alpha_i^{\vee}\}_{i\in I})$ associated to the extended Cartan matrix $A$, for $I=\{0,\dots,\ell\}$. 
If $\mathring{G}$ is simply-connected, the associated Kac--Moody group $G_{\mathcal D}$ is the loop group of $\mathring{G}$ with $k$-points
$$G_{\mathcal D}(k)=\mathring{G}(k((t,t^{-1}))).$$

To obtain a free root datum, let $X=\mathring X\oplus \Z \delta$ and $X^\vee=\mathring X^\vee \oplus \Z d$, where $\langle \delta, d\rangle=1$ and $\delta$ and $d$ are zero on $\mathring X^\vee$ and $\mathring X$, respectively. Moreover, modify the definition of the root $\alpha_0$ in $\mathcal D$ to $\alpha_0=-\theta+\delta\in X.$ The Kac--Moody root datum
$\mathcal D_{\text{free}}=(X, \{\alpha_i\}_{i\in I}, \{\alpha_i^{\vee}\}_{i\in I})$ is free. If $\mathring{G}$ is simply-connected, the associated Kac--Moody group $G_{\mathcal D_{\text{free}}}$ is the extension of the loop group by loop rotations with $k$-points
$$G_{\mathcal D_{\text{free}}}(k)=\mathring{G}(k((t,t^{-1})))\rtimes \mathbb{G}_m(k).$$
Here, $z\in \mathbb{G}_m(k)$ acts on $P(t)\in k((t,t^{-1}))$ via $z\cdot P(t)=P(zt)$.
If $\mathring{\mathcal{D}}$ is of simply-connected type then so is $\mathcal D_{\text{free}}.$

Dually, to obtain a cofree root datum, let $X=\mathring X\oplus \Z K^*$ and $X^\vee=\mathring X^\vee \oplus \Z K$, where $\langle K^*, K\rangle=1$ and $K^*$ and $K$ are zero on $\mathring X^\vee$ and $\mathring X$, respectively. Moreover, modify the definition of the coroot $\alpha_0^\vee$ in $\mathcal D$ to $\alpha_0^\vee=-\theta^\vee+K\in X^\vee.$ The Kac--Moody root datum
$\mathcal D_{\text{cofree}}=(X, \{\alpha_i\}_{i\in I}, \{\alpha_i^{\vee}\}_{i\in I})$ is cofree. If $\mathring{G}$ is simply-connected, the associated Kac--Moody group $G_{\mathcal D_{\text{cofree}}}$ is a central extension of the loop group of $\mathring G.$ If $\mathring{\mathcal{D}}$ is of adjoint type then so is $\mathcal D_{\text{cofree}}.$
\subsection{$K$-theory Soergel bimodules}\label{subsectionSoergelBimodules}
 We now explain how to attach a category of so-called $K$-theory Soergel bimodules to a Kac--Moody root datum $\mathcal{D}.$ Denote by $R=\Z[X]$ the \emph{character ring}. We use the exponential notation $e^\lambda\in R$ for $\lambda \in X.$
The ring $R$ inherits a natural action of the Weyl group $W=W(\mathcal{D}).$
For a simple reflection $s\in W$ we consider the $R$-bimodule
$$\BS(s)=R\otimes_{R^{s}}R$$
where $R^{s}$ are the $s$-invariants in $R.$
\begin{definition}
    The category of \emph{K-theory Bott--Samelson bimodules}, denoted by $\BSBimK_{\Dd}$, is the monoidal and additive subcategory of the category of $R$-bimodules generated by the objects $\BS(s_i)$ for $i\in I.$ The category of \emph{K-theory Soergel bimodules}, denoted by $\SBimK_{\Dd}$, is the idempotent-closed subcategory generated by $\BSBimK_{\Dd}.$
\end{definition}
We also record the following useful fact.
For $w\in W$, we denote by $R_w$ the twisted $R$-bimodule with $r\cdot m\cdot r'=rw(r')m$ for $r,r'\in R$ and $m\in R_w.$
\begin{lemma}\label{lem:hombetweenstandardbimodules} For $w,w'\in W$ we have $R_w\otimes_R R_{w'}=R_{ww'}$. If $\mathcal{D}$ is free, then in the abelian category of $R$-bimodules we have
$$
\Hom_{R\otimes R}(R_{w'}, R_{w}) = \begin{cases}
    R & \text{if } w = w', \\
    0 & \text{if } w \neq w'.
\end{cases}
$$
\end{lemma}
\begin{proof}
The first statement is clear. For the second, it suffices to check the case $w'=e\neq w.$ Let $\phi\in\Hom_{R\otimes R}(R, R_{w})$ and $m=\phi(1)=\sum a_\lambda e^\lambda$. Then $rm=w(r)m$ for all $r\in R$. Let $\mu\in X$ with $w(\mu)\neq \mu$. Let $\nu=w(\mu)-\mu$. Then $e^{\nu}m=m$. Assume that $a_\lambda\neq 0$ for some $\lambda\in X.$ Then also $a_{\lambda+i\nu}\neq 0$ for all $i\in\Z$. This is a contradiction, since almost all coefficients of $m$ are zero. So $m=0.$
\end{proof}
We now collect some helpful properties on $K$-theory Soergel bimodules. For this, let $i\in I$, $s=s_i$, $\alpha=\alpha_i$ and $\alpha^\vee= \alpha^\vee_i.$ The \emph{Demazure operator}, see \cite{demazureNouvelleFormuleCaracteres1975}, associated to $s$ is defined by 
$$D_s\colon R\to R,\, e^\lambda\mapsto \frac{e^\lambda-e^{s( \lambda)-\alpha}}{1-e^{-\alpha}}.
$$
The operator is $R^s$-linear and restricts to the identity on $R^s$.
We also consider the operator
$$D^-_s=1-D_s\colon R\to R,\, e^\lambda\mapsto \frac{e^{\lambda}-e^{s(\lambda)}}{1-e^{\alpha}}.$$
We obtain the following special case of a theorem of Pittie--Steinberg, see \cite{steinbergTheoremPittie1975}.
\begin{lemma}
Assume that there is a $\varpi\in X$ such that $\langle \varpi, \alpha^\vee\rangle=1.$ Then as an $R^s$-module, $R=R^s\oplus e^\varpi R^s$ with splitting given by the isomorphism
$$(D_s, D^-_s): R\to R^s\oplus e^{-\varpi} R^s.$$
This yields a splitting of the Bott--Samelson module
$$\BS(s)=R\otimes_{R^s}R=(R\otimes 1)\oplus (R\otimes e^{-\varpi})$$
as a left $R$-module.
\end{lemma}
We now consider the map of left $R$-modules
$$\op{gr}=(\op{gr}_e,\op{gr}_s): \BS(s)=R\otimes_{R^s}R\to R_e\oplus R_s,\, r\otimes r'\mapsto (rr', rs(r')).$$
\begin{lemma}\label{lemInjectivity}
 Assume that there is a $\varpi\in X$ such that $\langle \varpi, \alpha^\vee\rangle=1.$ Then, $\op{gr}$ is injective and its image, as a left $R$-module, is
$$\op{im}(\op{gr})=R(1,1)\oplus R(e^{-\varpi},e^{-s\varpi}).$$ This implies using $\varpi-s\varpi=\alpha$ that there is a map of short exact sequences of left $R$-modules
\[\begin{tikzcd}
	{R(1-e^\alpha)} & {\op{im}(\op{gr})} & {R_e} \\
	{R_s} & {R_e\oplus R_s} & {R_e.}
	\arrow[from=1-1, to=1-2]
	\arrow[hook, from=1-1, to=2-1]
	\arrow[from=1-2, to=1-3]
	\arrow[hook, from=1-2, to=2-2]
	\arrow[Rightarrow, no head, from=1-3, to=2-3]
	\arrow[from=2-1, to=2-2]
	\arrow[from=2-2, to=2-3]
\end{tikzcd}\]
\end{lemma}
We will also make use of the following fact.
\begin{lemma}\label{lemRingReduced}
Let $0\neq \lambda\in X$, then the ring $R/(1 - e^{\lambda})$ is reduced.
\end{lemma}
\begin{proof}
    We may choose an isomorphism $X\cong\mathbb{Z}v\oplus \Z^k$ such that $\lambda=nv$ for some $n\neq 0$. Then 
    $R/(1 - e^{\lambda})\cong \mathbb{Z}[\mathbb{Z}/n\oplus \mathbb{Z}^k]$ is reduced since it is the integral group ring of a commutative group, see \cite[Proposition 2.2]{mayGroupAlgebrasFinitely1976}.
\end{proof}

\section{Categories of $K$-motives}\label{sec:kmotives}
This section records definitions and basic properties of (reduced) $K$-motives $\DKbig$ on linearly reductive (ind-pro-)stacks.
\subsection{Quotient stacks} Fix a base scheme $\point=\Spec(k)$
for an algebraically closed field $k$ of characteristic zero.  A stack $\Xx/k$ is \emph{\nice} if it can be represented as a quotient stack $X/G$ where
\begin{enumerate}
    \item $G/k$ is a linearly reductive group,
    \item $X/k$ is a $G$-quasi-projective scheme, that is, it admits an embedding in a projectivised representation of $G$, and is of finite type.
\end{enumerate}
This definition includes quotients $X/G$ for linear algebraic groups $G\subset \op{GL}_n$ which are not linearly reductive, since $X/G\cong X\times^G\op{GL}_n/\op{GL}_n$ and $\op{GL}_n$ is linearly reductive.

We denote by $\Bc G=\point\!/G$ the classifying stack of a group $G.$
\subsection{$K$-motives on stacks}
For a \nice\ stack $\Xx$ we denote by  $\SH(\Xx)$ the stable motivic homotopy category associated to $\Xx$ as defined in \cite{hoyoisSixOperationsEquivariant2017}. 

By \cite{hoyoisCdhDescentEquivariant2020}, there is a ring spectrum $\KGL_\Xx$ in $\SH(\Xx)$ representing $\A^1$-homotopy invariant algebraic $K$-theory. The category of $K$-motives on $\Xx$ is defined by 
$$\DKbig(\Xx)=\Mod_{\KGL_\Xx}(\SH(\Xx)).$$
In \cite{hoyoisSixOperationsEquivariant2017,khanGeneralizedCohomologyTheories2022} it is shown that there is a functor between the category of correspondences of linearly reductive stacks with quasi-projective morphisms
\begin{align*}
    \op{Corr}(\op{Stk})&\to \op{Pr}: \Xx \mapsto \DKbig(\Xx),\, (\Xx\xleftarrow{f} \Zz \xrightarrow{g} \Yy)\mapsto g_!f^*.
\end{align*}
So there are the following functors
\begin{enumerate}
    \item $f^*,f_*,f_!,f^!$ for quasi-projective morphisms,
    \item and bifunctors $\iHom,\otimes.$ 
\end{enumerate}
These fulfill the usual axioms of a six functor formalism, such as base change, localization and projection formulae. We note that for $f$ smooth, there is an equivalence $f^*\simeq f^!$. This is related to \emph{Bott periodicity} for $K$-theory.

The category of $K$-motives computes algebraic $K$-theory and $G$-theory. For a stack $\Xx$ denote by 
$K(\Xx)$ and $G(\Xx)$ the $K$-theory spectra of the category of perfect complexes and coherent sheaves, respectively. If $\Xx/k$ is a smooth \nice\ stack and $f:\Yy\to\Xx$ is a quasi-projective map of \nice\ stacks, there are equivalences of spectra (in the first case of ring spectra)
\begin{align}
    \Map_{\DKbig(\Xx)}(\un,\un) &\simeq \op{K}(\Xx) \text{ and}\label{eq:DKcomputesKtheory}\\
    \Map_{\DKbig(\Yy)}(\un,f^!\un)& \simeq \op{G}(\Yy)\label{eq:DKcomputesGtheory},
\end{align}
see \cite[below Definition 5.1]{hoyoisCdhDescentEquivariant2020} and
 \cite[Remark 5.7]{hoyoisCdhDescentEquivariant2020}. 

\subsection{Reduced $K$-motives} \label{sec:reducedkmotives}
Since we are concerned with the Hecke category of a split reductive/Kac--Moody group---which is already defined over the integers---the constructions carried out here should be insensitive to the choice of base scheme $\point=\Spec(k).$ In particular, the higher algebraic $K$-theory $K_{>0}(\point)$ is irrelevant for us. 

To remove the higher $K$-groups, we pass to the category of \emph{reduced $K$-motives} $$\DK(\Xx)=\DKbig(\Xx)\otimes_{\op{K}(\point)}\op{K}_{0}(\point)$$
which is the Lurie tensor product of $\DKbig(\Xx)$ over the $K$-theory spectrum $\op{K}(\point)$ with $\op{K}_0(\point)=\Z.$ Here, the action of  $\op{K}(\point)$ on $\DKbig(\Xx)$ arises in the following way. We use that $\Xx$ is a linear reductive stack and hence of the form $X/G$ and denote by $f:X/G\to \Bc G$ the projection, which is quasi-projective by assumption. Then we obtain a pullback map $\op{K}(\point)\to \op{K}(\Bc G)=\Map_{\DKbig(\Bc G)}(\un,\un)\stackrel{f^*}{\to} \Map_{\DKbig(\Xx)}(\un,\un)$. The latter ring spectrum naturally acts on $\DKbig(\Xx)$.

The precise definition is explained in \cite{eberhardtIntegralMotivicSheaves2023} in the context of motivic sheaves $\DM$ and in \cite[Section 2.3]{eberhardtKmotivesSpringerTheory2024} for $K$-motives $\DKbig$. There it is shown that $\DK$ inherits the six functor formalism from $\DKbig.$ Furthermore, $\DK$ computes \emph{reduced $K$-theory} and \emph{reduced $G$-theory}
\begin{align}
    \Map_{\DK(\Xx)}(\un,\un) &\simeq \op{K}(\Xx)_\red = \op{K}(\Xx)\otimes_{\op{K}(\point)}\op{K}_{0}(\point) \text{ and}\label{eq:DKredcomputesredKtheory}\\
    \Map_{\DK(\Yy)}(\un,f^!\un)& \simeq \op{G}(\Yy)_\red = \op{G}(\Yy)\otimes_{\op{K}(\point)}\op{K}_{0}(\point)\label{eq:DKredcomputesredGtheory},
\end{align}
for $\Xx/k$ smooth and $f:\Yy\to\Xx$ a quasi-projective map of \nice\ stacks. Reduced $G$-theory is well-behaved for sufficiently `cellular' stacks.
\begin{definition}\cite[Definition 3.4]{eberhardtKmotivesSpringerTheory2024}
A stack $\Zz/k$ fulfills property $(C)$ if it admits a filtration
$\Zz=\Zz^n\supset \Zz^{n-1}\supset\cdots\Zz^1\supset \Zz^0=\emptyset$
into closed substacks such that 
$\Vv^i=\Zz^i-\Zz^{i-1}$ is a vector bundle over some classifying stack 
$\Bc(G_i\ltimes U_i)$ for $G_i$ reductive and $U_i$ unipotent.
\end{definition}
\begin{lemma}\label{lem:propertyCimpliesR} \cite[Lemma 3.3]{eberhardtKmotivesSpringerTheory2024} For a stack $\Zz/k$ fulfilling property $(C)$ there is an equivalence $\op{G}_{0}(\Zz)\stackrel{\sim}{\to}\op{G}(\Zz)_\red.$ If $\Zz$ is moreover smooth, $\op{K}$-theory and $\op{G}$-theory agree and $\op{K}_{0}(\Zz)\stackrel{\sim}{\to}\op{K}(\Zz)_\red.$
\end{lemma}
The categories of (reduced) $K$-motives $\DKbig(\Xx)$ and $\DK(\Xx)$ are too large for our purposes, since they for example contain the information about the (reduced) $G$-theory of $\Yy$ for all quasi-projective maps $f:\Yy\to \Xx$. For this reason, we will often restrict our attention to the stable subcategories generated by the tensor unit $\un$ which we think of as categories of locally constant objects.
\begin{example}\label{ex:unipotentgroupandtorus}
Let us consider two important special cases.
\begin{enumerate}
\item For $\Xx=\Bc U$, where $U$ is a unipotent group, the stable subcategory generated by the constant object $$\DK(\Bc U)\supset \langle \un \rangle_{\operatorname{stb}}\stackrel{\sim}{\to}\D^b(\Z)$$ is equivalent to the bounded derived category of $\Z$-modules. This follows from the fact that $K_0(\Bc U) \cong \Z$.
\item  For $\Xx=\Bc T$, where $T$ is a torus, the stable subcategory generated by the constant object 
$$\DK(\Bc T)\supset \langle \un \rangle_{\operatorname{stb}}\stackrel{\sim}{\to}\D^b(\Z[X(T)])$$ is equivalent to the bounded derived category of $\Z[X(T)]$-modules, where $X(T)$ is the character lattice of $T$. Here, we use that $K_0(\Bc T)=\Z[X(T)].$
\end{enumerate}
\end{example}

\subsection{Extension to certain ind-pro-stacks}\label{sec:extensiontoindstacks}
We now explain how to extend the formalism to certain ind-pro-stacks, keeping the Hecke stack $B\bs G/B$ of a Kac--Moody group in mind.
\begin{definition}
   A pro-\nice pro-stack $\Xx=\lim_{m\geq 0} \Xx_m$ is a pro-stack such that each $\Xx_m$ is \nice. We call a pro-\nice pro-stack $\Xx$ \niceprostack, if the transition maps $p_m\colon \Xx_{m}\to \Xx_{m-1}$ are smooth and quasi-projective.
\end{definition}
We can extend the formalism of $K$-motives to these stacks in the following way.
\begin{definition}
    The category of (reduced) $K$-motives on a \niceprostack pro-stack $\Xx=\lim_{m\geq 0} \Xx_m$ is defined as the limit
    $$\DKbig_{(\red)}(\Xx)=\lim_{m}{}_{\hspace{-0.1em}*}\DKbig_{(\red)}(\Xx_m)$$
    along the $*$-pushforwards of the transition maps.
\end{definition}
\begin{remark}\label{rem:independenceofpresentations}
    A pro-quasi-projective morphism $f$ from a \niceprostack pro-stack $\Xx=\lim_{m\geq 0} \Xx_m$ to a \nice stack $\Yy$ factors over some $\Xx_i\to \Yy$. This implies that the definition of $\DKbig_{(\red)}(\Xx)$ does not depend on the presentation of $\Xx$ in the category of stacks with quasi-projective maps.
\end{remark}
\begin{lemma}\label{lem:defprokmotives}
    For a \niceprostack pro-stack $\Xx=\lim_{m\geq 0} \Xx_m$, there is an equivalence
    $$\DKbig_{(\red)}(\Xx)=\lim_{m}{}_{\hspace{-0.1em}*}\DKbig_{(\red)}(\Xx_m)\simeq \underset{m}{\colim}^!\DKbig_{(\red)}(\Xx_m)$$
    with the colimit in $\mathrm{Pr}$ along the $!$-pullback maps.
    Moreover, $\DKbig_{(\red)}(\Xx)$ is compactly generated and inherits the functors $\otimes$, $\iHom$ and for pro-quasi-projective maps $f:\Xx\to \Yy$ the functors $f_*$, $f^!$ and $f^*$.
\end{lemma}
\begin{proof}
    The first statement follows since there is an adjunction
    $$p_m^*=p_m^!\colon \DKbig_{(\red)}(\Xx_{m-1})\leftrightarrows  \DKbig_{(\red)}(\Xx_{m}): p_{m,*}$$
    and the comparison of limits and colimit, see \cite[5.5.3.3]{lurieHigherToposTheory2009}. The existence of the functors follows from the compatibility with the (co)-limits.
\end{proof}
To also define a $!$-pushforward functor, we need the following restriction.
\begin{definition}
Consider pro-stacks $\Xx=\lim_{m\geq 0} \Xx_m$ and $\Yy=\lim_{n\geq 0} \Yy_n$. We call a morphism of pro-stacks $f:\Xx\to \Yy$ \emph{Cartesian} if for each $m$, there are indices $m'\geq m''$ such that $f$ factors through a Cartesian diagram
\[\begin{tikzcd}
	{\Xx_m} & {\Xx_{m-1}} \\
	{\Yy_{m'}} & {\Yy_{m''}}.
	\arrow[from=1-1, to=1-2]
	\arrow[from=1-1, to=2-1]
	\arrow[from=1-2, to=2-2]
	\arrow[from=2-1, to=2-2]
\end{tikzcd}\]
\end{definition}
Using base change, we obtain the following.
\begin{lemma}\label{lem:existenceflowershriekpro}
For a Cartesian morphism of \niceprostack pro-stacks $f:\Xx\to \Yy$,  the functor $f^!$ has a left adjoint $f_!$  induced by the functors $f_!:\DKbig_{(\red)}(\Xx_n)\to \DKbig_{(\red)}(\Yy_m)$ where $m\geq n$ such that $f|_{\Xx_n}$ factors through $\Yy_m$.
\end{lemma}
We will also use the following statement that is relevant for quotients by pro-unipotent groups.
\begin{lemma}\label{lem:fullyfaithfulpro}
Let  $\Xx=\lim_{m\geq 0} \Xx_m$ be a \niceprostack pro-stack where we assume additionally that the transition maps $p_m$ are torsors over vector bundles. In this case $p_m^!=p_m^*$ is fully faithful and hence the natural insertion map 
$$\op{ins}_0: \DKbig_{(\red)}(\Xx_0)\to \DKbig_{(\red)}(\Xx)$$
is fully faithful as well.
\end{lemma}
Now we extend our formalism to ind-pro-stacks.
\begin{definition}
    Let $\Xx=\colim \Xx_n$ be an ind-linearly-reductive ind-pro-stack. We call $\Xx$ \niceindprostack, if the transition maps $\iota_n: \Xx_n\to \Xx_{n+1}$ are Cartesian pro-proper.
\end{definition}
\begin{definition}
    The category of (reduced) $K$-motives on a \niceindprostack ind-pro-stack $\Xx=\colim \Xx_n$ is defined as the colimit in $\mathrm{Pr}$
    $$\DKbig_{(\red)}(\Xx)=\underset{m}{\colim}{}_{\hspace{0.05em}*}\DKbig_{(\red)}(\Xx_m)$$
    along the $*$-pushforward maps.
\end{definition}
\begin{remark}
    As in \Cref{rem:independenceofpresentations}, this definition does not depend on the presentation of $\Xx$ in the category of pro-stacks with pro-quasi-projective morphisms.
\end{remark}
\begin{lemma}
For a \niceindprostack ind-pro-stack $\Xx=\colim \Xx_n$, there is an equivalence of categories
$$\DKbig_{(\red)}(\Xx)=\underset{m}{\colim}{}_{\hspace{0.05em}*}\DKbig_{(\red)}(\Xx_m)\simeq \lim_{m}{}^{\hspace{-0.1em}!}\DKbig_{(\red)}(\Xx_m)$$
with the limit along the $!$-pullbacks. Moreover, $\DKbig_{(\red)}(\Xx)$ inherits the functors $\otimes$, $\iHom$ and for ind-pro-quasi-projective morphisms $f:\Xx\to \Yy$ the functors $f_*$, $f^!.$ 
\end{lemma}
\begin{proof}
    This works as in Lemma \ref{lem:defprokmotives} using that there is an adjunction
    \begin{align*}\iota_{n,*}&=\iota_{n,!}\colon \DKbig_{(\red)}(\Xx_{n})\leftrightarrows  \DKbig_{(\red)}(\Xx_{n+1}): \iota_{n}^!\qedhere\end{align*}
\end{proof}
\begin{definition}
We call a morphism of ind-pro-stacks $f:\Xx\to \Yy$ Cartesian if for all $n$ there are $n'\geq n''$ such $f$ factors through a Cartesian diagram of pro-stacks
\[\begin{tikzcd}
	{\Xx_n} & {\Xx_{n+1}} \\
	{\Yy_{n'}} & {\Yy_{n''}}
	\arrow[from=1-1, to=1-2]
	\arrow[from=1-1, to=2-1]
	\arrow[from=1-2, to=2-2]
	\arrow[from=2-1, to=2-2]
\end{tikzcd}\]
\end{definition}
\begin{lemma}
For an ind-pro-quasi-projective morphism $f:\Xx\to \Yy$ of \niceindprostack ind-pro-stacks, the functor $f_*$ has a left adjoint $f^*$ induced by the functors $f^*:\DKbig_{(\red)}(\Yy_{n'})\to \DKbig_{(\red)}(\Xx_n)$. Similarly, the functor $f_!$ as considered in Lemma \ref{lem:existenceflowershriekpro} extends to \niceindprostack ind-pro-stacks, if $f$ is ind-Cartesian.
\end{lemma}
Similarly to \Cref{lem:fullyfaithfulpro} we obtain the following.
\begin{lemma}
    If in a \niceindprostack ind-pro-stack $\Xx=\colim \Xx_n$ the transition maps $\iota_n$ are pro-closed immersions, $\iota_{n,*}$ is fully faithful and hence the insertion maps
    $$\op{ins}_n: \DKbig_{(\red)}(\Xx_n)\to \DKbig_{(\red)}(\Xx)$$
    are fully faithful.
\end{lemma}
In summary, let 
$\Xx=\colim_n \lim_m \Xx_{n,m}$ be an ind-pro stack with a presentation of the form
\begin{enumerate}
    \item $\Xx_{n,m}$ is a \nice stack,
    \item $p_{n,m}: \Xx_{n,m}\to \Xx_{n,m-1}$ is a torsor under a vector bundle and
    \item $\iota_{n,m}:\Xx_{n,m}\to \Xx_{n+1,\varphi(n,m)}$ is a closed immersion.
\end{enumerate}
Then we can consider the category of reduced $K$-motives $\DKbig_{(\red)}(\Xx)$ on $\Xx$ and for each $n$ there are fully faithful insertion functors
$$\op{ins}_n\op{ins}_0: \DKbig_{(\red)}(\Xx_{0,n})\to \DKbig_{(\red)}(\Xx).$$
For most purposes, we can hence restrict ourselves to considering $K$-motives on the \nice stacks $\Xx_{0,n}.$

\section{Soergel description of $K$-motives on flag varieties} \label{sec:soergelktheoretichecke}
In this section we will give a Soergel-theoretic description of the $K$-theoretic Hecke category $\HeckeK_{\Dd}$ which consists of $K$-motives on the Hecke stack $B\backslash G/B$ of a Kac--Moody group $G$ associated to a Kac--Moody root datum $\Dd$.
\subsection{$K$-motives on the Hecke stack}
Recall that $k$ is an algebraically closed field of characteristic $0$ and $\point=\Spec(k).$ Denote by $G\supset B\supset T$ the Kac--Moody group together with a Borel subgroup and maximal torus associated to the Kac-Moody datum $\mathcal{D}$, see \Cref{sec:kacmoodyrootdatum} and \Cref{sec:kacmoodygroups}. Denote by $W\supset S$ the Weyl group with the set of simple reflections.

The \emph{Hecke stack} $B\bs G/B$ is ind-pro-\nice, see \Cref{sec:kacmoodyflagvarieties}. We can hence consider the category of reduced $K$-motives $\DK(B\bs G/B)$, see \Cref{sec:extensiontoindstacks}. 
For $w\in W$, we denote by $i_{w}: B\bs BwB/B\to B\bs G/B$ the inclusion. We denote by $$\Delta_w=i_{w,!}\un, \nabla_w=i_{w,*}\un\in \DK(B\bs G/B)$$ the corresponding \emph{standard} and \emph{costandard objects}. 
\begin{definition} The \emph{$K$-theoretic Hecke category} is the stable subcategory generated by the standard objects
$$\Hecke=\HeckeK_{\Dd}=\genbuild{\Delta_w}{w\in W}_{\op{stb}}\subset \DK(B\bs G/B).$$
\end{definition}
\begin{remark}
    \begin{enumerate}
        \item The map $T\bs G/B\to B\bs G/B$ is a pro-affine bundle. This implies that the pullback $\DK(B\bs G/B)\to \DK(T\bs G/B)$ is fully faithful with essential image objects which are locally constant along $B$-orbits on  $T\bs G/B$. While this perspective is sometimes helpful since $T\bs G/B$ is simply an ind-\nice stack, the asymmetry is problematic when defining the monoidal structure.
        \item We note that the $K$-theoretic Hecke category $\HeckeK_{\Dd}$ is much smaller than the category of all reduced $K$-motives $\DK(B\bs G/B)$. The former only contains objects that are `locally constant' along the Bruhat cells.
        \item If $G=B=T$ is a torus, then there is a single standard object and $B\backslash G/B=\point/T.$ By \Cref{ex:unipotentgroupandtorus} there is an equivalence $$\HeckeK_{\Dd}\stackrel{\sim}{\to}\D^b(\Z[X(T)]).$$
    \end{enumerate}
\end{remark}

\subsection{Monoidal structure} The Hecke stack is a monoid in the category of correspondences of ind-pro-\nice stacks with multiplication given by the convolution diagram
\begin{equation}\label{eq:convolutiondiagram}
    \begin{tikzcd}
        & {B \bs G \times^B G/B} & {B \bs G/B} \\
        {B \bs G/B} && {B \bs G/B}
        \arrow["m", from=1-2, to=1-3]
        \arrow["{p_2}"', from=1-2, to=2-3]
        \arrow["{p_1}", from=1-2, to=2-1]
    \end{tikzcd}    
\end{equation}
The maps $p_1,p_2$ and $m$ are ind-pro-proper with fiber $G/B.$ The inclusion map $i_e: B\bs B /B\to B\bs G/B$ is the unit. This turns $\DK(B\backslash G/B)$ into a monoidal category: There are two monoidal products defined by
\begin{align*}
    A\conv B &= m_*(p_1^*A\otimes p_2^*B)=m_!(p_1^*A\otimes p_2^*B)\text{ and }\\
    A\conv^! B &= m_*(p_1^!A\otimes p_2^!B)=m_!(p_1^!A\otimes p_2^!B)
\end{align*}
for $A,B\in \DK(B\backslash G/B)$. For both products, the unit object is $i_{e,*}\un$.

It will turn out that $\conv$ and $\conv^!$ are equivalent when restricted to the Hecke category $\Hecke$, see \Cref{thm:rigiditykmotivichecke}. This is immediate if $G$ is a reductive group, since then $G/B$ is smooth and $p_i^*\cong p_i^!$ for $i=1,2.$

We now show that the monoidal structures restrict to $\Hecke.$ For $w\in W$, denote 
$$\Hh_{\leq w}=\genbuild{\Delta_x}{x\in W, x\leq w}_{\op{stb}}\subset \Hecke.$$ 
For $s\in S,$  let $P_s=B\cup BsB$ be the corresponding standard parabolic and $k_s:B\bs P_s/B\to B\bs G/B$ the inclusion. Denote $\IK_s=k_{s,!}\un_s$ and $\IK_e=\Delta_e.$
\begin{lemma}\label{lem:convolutionsandstandards} Let $x,y\in W$ and $s\in S.$ Then
    \begin{enumerate}
        \item $\Delta_x\conv\Delta_y=\Delta_{xy}$ and $\nabla_x\conv^!\nabla_y=\nabla_{xy}$ if $\ell(xy)=\ell(x)+\ell(y),$
        \item $\genbuildexplicit{\IK_e,\IK_s}=\genbuildexplicit{\Delta_e,\Delta_s}=\genbuildexplicit{\nabla_e,\nabla_s}=\Hh_{\leq s},$
        \item $\IK_s\conv \IK_s=\IK_s\conv^! \IK_s\cong \IK_s \oplus \IK_s$ and
        \item $\Hh_{\leq s}\conv \Hh_{\leq s}=\Hh_{\leq s}\conv^! \Hh_{\leq s}\subset \Hh_{\leq s}.$
    \end{enumerate}
\end{lemma}
\begin{proof} (1) follows from base change and that multiplication yields an isomorphism
    $BxB\times^B ByB\stackrel{\sim}{\to} BxyB.$
    
    (2) follows from the localization cofiber sequences 
\[\begin{tikzcd}[column sep=scriptsize,row sep=tiny]
	{\Delta_e=E_e} & {E_s} & {\Delta_s} \\
	{\nabla_s} & {E_s} & {\nabla_e=\IK_e}
	\arrow[from=1-2, to=1-3]
	\arrow[from=1-1, to=1-2]
	\arrow[from=2-1, to=2-2]
	\arrow[from=2-2, to=2-3]
\end{tikzcd}\]
induced by the decomposition $P_s=B\cup BsB$.

(3) follows from the projective bundle formula applied to $m: B\bs P_s\times^B P_s/B\to B\bs P_s/B$ which has fiber $P_s/B\cong \mathbb{P}^1.$

(4) follows from (2), (3) and the fact that $\Delta_e=\IK_e$ is the unit for convolution.
\end{proof}
\begin{corollary}
    The category $\Hecke$ is stable under $\conv$ and $\conv^!.$
\end{corollary}
\begin{proof}
By induction it suffices to show that $\Delta_w\conv \Delta_s\in \Hecke$ for $w\in W$ and $s\in S$. If $ws>x$, then $\Delta_w\conv \Delta_s=\Delta_{ws}\in \Hecke$. If $ws<w$, $\Delta_w\conv \Delta_s=\Delta_{ws}\conv \Delta_s\conv \Delta_s\in \langle \Delta_{ws}, \Delta_{w}\rangle\subset \Hecke.$ The statement about $\conv^!$ follows similarly.
\end{proof}
We will also need the following statement for standard parabolic subgroups.
\begin{lemma}\label{lem:convolutionfiniteparabolic}
Let $B\subset P\subset G$ be a parabolic subgroup with a Kac--Moody datum of finite type. Let $w_P\in W$ the longest element in the Weyl group $W_P$ of $P$. Then the two convolution products $\conv$ and $\conv^!$ become equivalent when restricted to $\Hecke_{\leq w_P}.$
\end{lemma}
\begin{proof}
    This follows since the monoidal structure on $B\bs G/B$ restricts to $B\bs P/B$ and since $P/B$ as the flag variety of a reductive group is smooth.
\end{proof}

\subsection{Pure Objects} We define the category of \emph{pure} objects in $\Hecke$ as the additive, idempotent-closed monoidal (with respect to $\conv$) subcategory 
$$\Heckepure=\genbuild{\IK_s}{s\in S}_{\conv, \inplus,\cong}\subset \Hecke$$
generated by the objects $\IK_s.$ 
For $\underline{x}=(s_1,\dots,s_n)\in S^n$ we abbreviate the associated \emph{Bott--Samelson $K$-motive} by
$$E_{\underline{x}}=E_{s_1}\conv\cdots\conv E_{s_n}.$$
\begin{lemma}\label{lem:bottsamelson} For a sequence $\underline{x}=(s_1,\dots,s_n)\in S^n,$ denote 
    $$\pi_{\underline{x}}:\BS_{\underline{x}}=B\bs P_{s_1}\times^B\cdots\times^BP_{s_n}/B\to B\bs G/B.$$
Then there are isomorphisms
    $$\IK_{\underline{x}}=\IK_{s_1}\conv\cdots\conv \IK_{s_n}\cong\pi_{\underline{x},!}\un=\pi_{\underline{x},*}\un\cong \IK_{s_1}\conv^!\cdots\conv^! \IK_{s_n}.$$
\end{lemma}
\begin{proof}
This follows similarly to \cite[Section 3.2]{soergelRelationIntersectionCohomology2000a}.
\end{proof}
\begin{remark}
    \Cref{lem:bottsamelson} justifies the name \emph{pure} for the objects in $\Heckepure$: In the yoga of weights, pure objects correspond to (summands) of smooth and projective varieties.
\end{remark}
\begin{corollary}\label{cor:puregeneration} The category of pure objects is generated by the Bott--Samelson objects $\IK_{\underline{x}}$ as an additive and idempotent closed subcategory and one may replace $\conv$ by $\conv^!$ in its definition $$\Heckepure=\genbuild{\IK_{\underline{x}}}{\underline{x}\in S^n}_{\inplus,\cong}=\genbuild{\IK_s}{s\in S}_{\conv^!, \inplus,\cong}.$$
\end{corollary}
\begin{corollary}\label{cor:puregeneratehecke}
Pure objects stably generate the Hecke category:
$$\genbuildexplicit{\Heckepure}_{\op{stb}}=\Hecke.$$
\end{corollary}
\subsection{Pointwise Purity and Formality}
The pure objects in the Hecke category fulfill an additional pointwise purity property which allows to show a formality result. The following definitions and statements are a variation of \cite[Section 6]{soergelPerverseMotivesGraded2018}.
\begin{definition}
    For $?\in\{!,*\}$, we call an object $M\in \Hecke$ pointwise $?$-pure if for all $w\in W$ the restrictions $i^?_wM\in \DK(B\bs BwB/B)$ are a finite direct sum of the constant object $\un.$ An object is simply called pointwise pure if it is both pointwise $*$-pure and $!$-pure.
\end{definition}
There is the following equivalent description of pointwise pure objects.
\begin{lemma}\label{rem:pointwisepureyieldsfiltrations}
    Pointwise $!$-pure objects $E\in \Hecke$ are exactly the objects that admit a filtration whose associated graded are costandard objects, that is, there is a sequence of objects $E^k\in \Hecke$ for $k=0,\dots,n$ such that $E^0=0$, $E^n=E$ and there are cofiber sequences of the form
    \[\begin{tikzcd}[column sep=scriptsize,row sep=tiny]
	{E^{k-1}} & {E^k} & {\nabla_{w_k}}
	\arrow[from=1-2, to=1-3]
	\arrow[from=1-1, to=1-2]
\end{tikzcd}\]
    for $w_k\in W$ and $k=1,\dots n.$
    Dually, pointwise $*$-pure objects are precisely the objects that admit a cofiltration whose associated graded are standard objects.
\end{lemma}
\begin{proof}
Assume that $E\in \Hecke$ is pointwise $*$-pure. Choose a linear order on $W$ refining the Bruhat order. Then we obtain a cofiltration of $E$ by objects $E'^k=i_{\leq k,*}i_{\leq k}^!E$ where $i_{\leq k}$ is the closed embedding of the union of all Bruhat cells for $w_{k'}\in W$ with $k'\leq k.$ Then we obtain the localization cofiber sequence
\[\begin{tikzcd}[column sep=scriptsize,row sep=tiny]
	{E'^{k-1}} & {E'^k} & {i_{w_k,*}i_{w_k}^! E.}
	\arrow[from=1-2, to=1-3]
	\arrow[from=1-1, to=1-2]
\end{tikzcd}\]
By assumption, $i_{w_k}^! E$ is a finite direct sum of constant objects $\un$, so  $i_{w_k,*}i_{w_k}^! E$ is a finite direct sum of costandard objects $\nabla_{w_k}=i_{w_k,*}\un$. The desired filtration can hence be obtained by refinement.

Now assume that $E$ admits a filtration by objects $E^k$ whose associated graded are costandard objects. We show that $E$ is pointwise $!$-pure by induction on the length of the cofiltration. Let $w\in W$ and consider the cofiber sequence
\[\begin{tikzcd}[column sep=scriptsize,row sep=tiny]
	{i_w^!E^{k-1}} & {i_w^!E^k} & {i_w^!\Delta_{w_k}.}
	\arrow[from=1-2, to=1-3]
	\arrow[from=1-1, to=1-2]
\end{tikzcd}\]
By induction, $E^{k-1}$ is pointwise $!$-pure, so $i_w^!E^{k-1}$ is a finite direct sum of objects $\un$. If $w\neq w_k$, then $i_w^!\nabla_{w_k}=0$ and hence $i_w^!E^k=i_w^!E^{k-1}$ is also a finite direct sum of objects $\un$. If $w=w_k$, then $i_w^!\nabla_{w_k}=\un.$ Since
$$\Hom_{\DK(B\bs BwB/B)}(\un,\un[1])=\pi_{-1}(K(B\bs BwB/B)\otimes_{K(\point)}K_0(\point))=0,$$ the cofiber sequence splits and $i_w^!E^k=i_w^!E^{k-1}\oplus \un$ is also a finite direct sum of objects $\un$. Hence $E^k$ is pointwise $!$-pure.

The statement for pointwise $*$-pure objects follows by dual arguments. 
\end{proof}
An important consequence is the following formality of mapping spaces of pointwise pure objects.
\begin{lemma}\label{lem:pointwisepureformality}
    Let $M,N\in \Hecke$ be pointwise $*$-pure and pointwise $!$-pure, respectively. Then 
    $\Map(M,N)$ is concentrated in degree zero, so $\Map(M,N)= \Hom(M,N).$
\end{lemma}
\begin{proof}
 By base change and \eqref{eq:DKredcomputesredKtheory} we see that
    \begin{align*}
        \Map(\Delta_x, \nabla_y)&=K(B\bs (BxB\cap ByB)/B)\otimes_{K(\point)}K_0(\point)\\
        &=K_0(B\bs (BxB\cap ByB)/B)
    \end{align*}
    is concentrated in degree zero for all $x,y\in W.$
    The general statement follows by a double induction on the length of a cofiltration by standard modules of $M$ and filtration by costandard modules of $N$ (see \Cref{rem:pointwisepureyieldsfiltrations}) using the five lemma.
\end{proof}
We now show that pure objects in the Hecke category are also pointwise pure.
\begin{lemma}\label{lem:pureimpliespointwisepure}
    All objects in $\Heckepure$ are pointwise pure.
\end{lemma}
\begin{proof}
It suffices to show that the Bott-Samelson objects $E_{\underline{x}}=\pi_{\underline{x},!}\un=\pi_{\underline{x},*}\un$ (in the notation of \Cref{lem:bottsamelson}) are pointwise pure. We show that $E_{\underline{x}}$ is pointwise $!$-pure,  the pointwise $*$-purity follows by a dual argument. For this, we construct a filtration of $E_{\underline{x}}$ whose associated graded are costandard objects, see \Cref{rem:pointwisepureyieldsfiltrations}.

Similar to the proof of \Cref{rem:pointwisepureyieldsfiltrations} there is a filtration of $E_{\underline{x}}$ with associated graded of the form $i_{w,*}i_w^!\pi_{\underline{x},*}\un.$ The fiber $\pi_{\underline{x}}^{-1}(B\bs BwB /B)$ admits a stratification whose fibers are affine bundles over $B\bs BwB /B$, see \cite{hainesProofKazhdanLusztigPurity,haerterich2004tequivariantcohomologybottsamelsonvarieties}. Using this, we can refine the filtration such that the associated graded pieces are of the form $i_{w,*}p_*\un=\un$ where $p$ is an affine bundle over $B\bs BwB /B.$ Hence, the filtration has the desired form.
\end{proof}
Since pure objects stably generate the Hecke category, see \Cref{cor:puregeneratehecke}, with \Cref{prop:weightcomplexequivalence} and \Cref{lem:pointwisepureformality} we obtain the following formality result for the $K$-theoretic Hecke category.
\begin{corollary}\label{cor:formalityofkmotivichecke}
    There is an equivalence of monoidal categories between the $K$-theoretic Hecke category and the category of bounded chain complexes of pure objects in $\Hecke$
    $$\Hecke\stackrel{\sim}{\to}\Ch^b(\Heckepure).$$
\end{corollary}
\subsection{Duality} 
Our next goal is to show that the Hecke category $\Hecke$ is rigid.
For $M\in \DK(B\backslash G/B)$ we define
\begin{align*}
    \Dual(M)=\iHom(M,\omega) \text{ and }\DualI(M)=\Dual(\op{inv}^*M).
\end{align*}
where $\op{inv}$ is the inversion map of $G$ and $\omega=p^!(\un)$ and $p :B\bs G/B\to  B\bs \point/B$ is the projection.

We collect some properties of the duality functors.
\begin{lemma}\label{lem:propertiesduality}
Let $\underline{x}\in S^n$ and $\underline{x}^{\op{op}}$ the reversed sequence. Let $w\in W.$
    \begin{enumerate}
        \item There are isomorphism $\Dual(\IK_{\underline{x}})\cong \IK_{\underline{x}}$ and $\DualI(\IK_{\underline{x}})\cong \IK_{\underline{x}^{\op{op}}}$.
        \item $\Dual$ and $\DualI$ preserve $\Hh$ and $\Dual^2\cong \id\cong (\DualI)^2$ on $\Hh.$
        \item There is a natural map $\Dual(A)\conv^!\Dual(B)\to \Dual(A\conv B)$ that restricts to an isomorphism on $\Hh.$
        \item There is a natural map $\Dual(A\conv^! B)\to \Dual(A)\conv\Dual(B)$ that restricts to an isomorphism on $\Hh.$
        \item There is a natural map $\Dual^-(A\conv B)\to \Dual^-(B)\conv^!\Dual^-(A)$ that restricts to an isomorphism on $\Hh.$
    \end{enumerate}
\end{lemma}
\begin{proof}
    (1) follows from \Cref{lem:bottsamelson} as well as $\op{inv}^*(A\conv B)=\op{inv}^*(B)\conv \op{inv}^*(A)$ and $\op{inv}^*{E_s}=E_s$ for $s\in S.$

    (2) Follows from (1) using that $\Heckepure$ stably generates $\Hh.$
    
    (3) and (4) follow since $\Dual$ exchanges $!$ and $*$ and $\Dual^2=\id$ by (2). Similarly, (5) follows by keeping track of the extra $\op{inv}^*$.
\end{proof}

\begin{lemma}\label{lem:dualizability}
    Let $A,B,C\in\Hecke$, then there are natural isomorphisms
    $$\Hom(A\conv B, C)=\Hom(A,C\conv^! \DualI(B))=\Hom(B,\DualI(A)\conv^! C).$$
\end{lemma}
\begin{proof}
    We just show the first isomorphism. 
    
    For $A,B\in \Hecke$, using $\Dual^2=\id$, there is a natural isomorphism
    $$\Hom(A,B)\cong \Hom_{\DK(B\bs \point/B)}(p_!(A\otimes \Dual(B)),\un).$$
Using this and \Cref{lem:propertiesduality} we can hence write 
\begin{align*} 
    \Hom(A\conv B, C)&=\Hom_{\DK(B\bs \point/B)}(p_!(A\conv B\otimes \Dual(C)),\un)\text{ and }\\
    \Hom(A,C\conv^! \DualI(B))&=\Hom_{\DK(B\bs \point/B)}(p_!( A\otimes \Dual(C\conv^! \DualI(B))),\un)\\&=\Hom_{\DK(B\bs \point/B)}(p_!(A\otimes \Dual(C)\conv \op{inv}^*(B)),\un).
\end{align*}
Replacing $\Dual(C)$ by $C$, we hence need to show that there is a natural isomorphism
$$p_!((A\conv B)\otimes C)\cong p_!(A\otimes (C\conv \op{inv}^*(B))).$$
Using that $m_!=m_*$ as well as the projection formula, we get
\begin{align*}
    (A\conv B)\otimes C&=m_*(p_1^*A\otimes p_2^*B\otimes m^*C)\text{ and }\\
    A\otimes (C\conv \op{inv}^*(B))&=m_*(m^*A\otimes p_1^*C\otimes p_2^*\op{inv}^*B)
\end{align*}
Consider the isomorphism $r:B\bs G\times^B G/B, [x,y]\mapsto [xy,y^{-1}]$. Then $p_1r=m$, $p_2r=\op{inv}p_2$ and $mr=p_1$. Hence,
$$r^*(p_1^*A\otimes p_2^*B\otimes m^*C)=m^*A\otimes p_1^*C\otimes p_2^*\op{inv}^*B.$$
Using $r_*r^*=\id$, $pmr=pm$, $m_!=m_*$ and $r_!=r_*$, we obtain
\begin{align*}
    p_!((A\conv B)\otimes C)&=p_!m_*(p_1^*A\otimes p_2^*B\otimes m^*C)\\
    &=p_!m_*r_*r^*(p_1^*A\otimes p_2^*B\otimes m^*C)\\
    &=p_!m_*(m^*A\otimes p_1^*C\otimes p_2^*\op{inv}^*B)\\
    &=p_!(A\otimes (C\conv \op{inv}^*(B))).\qedhere
\end{align*}
\end{proof}
\begin{theorem}\label{thm:rigiditykmotivichecke}
    On $\Hecke$, there is a natural equivalence $\conv\simeq \conv^!$. Moreover, $\Hecke$ is rigid with left and right dual given by $\DualI$.
\end{theorem}
\begin{proof}
    To show this we will make use of \cite{boyarchenkoDualityFormalismSpirit2013}.
    In the notation of \cite{boyarchenkoDualityFormalismSpirit2013}, the object $\Delta_e$ is dualizing in the category $\Hecke$ with the monoidal structure given by the $*$-convolution. The duality functor is $\Dual^-$ using \Cref{lem:dualizability} where we use that $\Dual^-$ is an anti-equivalence by \Cref{lem:propertiesduality}. By \cite[Section 3.1]{boyarchenkoDualityFormalismSpirit2013} and using \Cref{lem:propertiesduality}, there is a natural map 
    $$X\conv Y\to \Dual^-(\Dual^-Y\conv \Dual^-X)=X\conv^!Y.$$ The objects $\IK_s$ are left and right dualizable by \Cref{lem:dualizability} using that on $B\bs P_s/B$ we have $\conv=\conv^!$ by \Cref{lem:convolutionfiniteparabolic}.  By \Cref{cor:puregeneration}, all objects in $\Heckepure$ are left and right dualizable as well. This implies by \cite[Lemma 3.4]{boyarchenkoDualityFormalismSpirit2013}, that the natural map between $*$ and $*^!$ is an isomorphism on $\Hecke$. The statement now follows from \cite[Corollary 3.6]{boyarchenkoDualityFormalismSpirit2013}.
\end{proof}
\subsection{Erweiterungssatz} 
We will now give a Soergel description for the pure objects in the Hecke category.
Let $R=\Z[X(T)].$ Then by \Cref{lem:propertyCimpliesR} we can identify
$$\Map_{\DK(B\bs \point /B)}(\un,\un)=K_0(T\bs \point /T)=R\otimes R.$$
Consider the functor
$$\Kyp: \DK(B\bs G/B)\to \D(R\otimes R),\, \Kyp(M)=\Map_{\DK(B\bs \point /B)}(\un,p_*M).$$
\begin{remark}
    The notation $\Kyp$ is inspired by the notation $\mathbb{H}$ for the hypercohomology functor which is used in the setting of (equivariant) constructible sheaves, see \cite{soergelKategorieMathcalPerverse1990,soergelCombinatoricsHarishChandraBimodules1992}.
\end{remark}

\begin{lemma}
    The functor $\Kyp$ is lax monoidal where $\DK(B\bs G/B)$ and $\D(R\otimes R)$ are equipped with the monoidal structures $\conv$ and $\otimes$, respectively.
\end{lemma}
\begin{proof}
    The functor $p_*$ is lax monoidal when we equip $\DK(B\bs \point /B)$ with the monoidal structure $\otimes$. To see this, we use that there is a natural map
    \begin{align*}
    p_*(-)\otimes p_*(-)&\to p_*p_{1,*}p_1^*(-)\otimes p_*p_{2,*}p_2^*(-)\\
    &=\overline{p}_*p_1^*(-)\otimes \overline{p}_*p_2^*(-)\\
    &\to \overline{p}_*(p_1^*(-)\otimes p_2^*(-))\\
    &= p_*m_*(p_1^*(-)\otimes p_2^*(-))\\
    &\stackrel{\sim}{\to} p_*m_!(p_1^*(-)\otimes p_2^*(-))=p_*(-\conv -)
    \end{align*}
    where we use the notation of \eqref{eq:convolutiondiagram} and denote $\overline{p}: B\bs G\times^B G/B\to B\bs \point /B.$ We then use that $\Map(\un,-)$ is lax monoidal.
\end{proof}
We now show that the functor $\Kyp$ evaluated at Bott--Samelson $K$-motives `computes' the $K$-theory of Bott--Samelson resolutions.
\begin{lemma}\label{lem:kypofbottsamelson}
 Let $\underline{x}\in S^n$, then there is an equivalence of $R$-bimodules
    $$\Kyp(\IK_{\underline{x}})\simeq K_0(\BS_{\underline{x}}).$$
\end{lemma}
\begin{proof}
This follows from \eqref{eq:DKredcomputesredKtheory}, \Cref{lem:bottsamelson} and \Cref{lem:propertyCimpliesR} using that $\BS_{\underline{x}}$ has property $(C)$ (to be precise a suitable pro-version of property $(C)$) and is (pro)-smooth.
\end{proof}
Under an appropriate assumption on the root datum, the functor $\Kyp$ is actually monoidal, when we equip $\D(R\otimes R)$ with the monoidal structure $-\otimes_R -$.
\begin{theorem}\label{thm:monoidalityofkyp}
    Assume that the Kac--Moody datum $\Dd$ is of simply-connected type. Let $\underline{x}\in S^n$ and $s\in S$.
    \begin{enumerate}
    \item The lax monoidal structure yields an isomorphism
    $$\Kyp(E_s\conv E_{\underline{x}})\stackrel{\sim}{\leftarrow}\Kyp(E_s)\otimes_R\Kyp(E_{\underline{x}})=R\otimes_{R^s}\Kyp(E_{\underline{x}}).$$
    \item When restricted to $\Hecke$, the functor $\Kyp$ is monoidal, where $\Hecke$ and $\D(R\otimes R)$ are equipped with the monoidal structures $\conv$ and $\otimes_R$, respectively.
    \item We have $\Kyp(\nabla_w)=R_w$ and $\Kyp(E_s)=R\otimes_{R^s}R$.
    \end{enumerate}
\end{theorem}
\begin{proof}
    (1) is shown in \cite[Section 4.4]{eberhardtKtheorySoergelBimodules2024}. Here we use that $\Dd$ is of simply-connected type.

    (2) is true  when restricting $\Kyp$ to the subcategory  $\Heckepure$ by induction and (1). Now use that $\Heckepure$ generates $\Hecke$ as a stable category.

    (3) is shown in \cite[Section 4.4]{eberhardtKtheorySoergelBimodules2024}.
\end{proof}
We are now ready to prove a $K$-theoretic version of Soergel's Erweiterungssatz \cite{soergelKategorieMathcalPerverse1990} for Kac--Moody groups. See \cite{eberhardtKtheorySoergelBimodules2024} for the case of reductive groups.
\begin{theorem}[Erweiterungssatz]\label{thm:erweiterungssatz} Assume that the Kac--Moody datum $\Dd$ is free and of simply-connected type. Then the functor $\Kyp: \Heckepure\to \D(R\otimes R)$ is fully faithful.
\end{theorem}
\begin{proof}
    Let $\underline{x}$ and $\underline{y}$ be sequences of simple reflections. It suffices to show that $\Kyp$ induces an isomorphism
    $$\Hom_{\Hecke}(\IK_{\underline{x}}, \IK_{\underline{y}})\to\Hom_{\D(R\otimes R)}(\Kyp(\IK_{\underline{y}}),\Kyp(\IK_{\underline{x}})).$$
    Since $\Kyp$ is monoidal by \Cref{thm:monoidalityofkyp} and the object $\IK_{\underline{y}}$ is left and right dualizable with dual $\IK_{\underline{y}^{op}}$ by \Cref{lem:propertiesduality} and \Cref{thm:rigiditykmotivichecke}, we obtain the commutative square

\[\begin{tikzcd}
	{\Hom_{\Hecke}(\IK_{\underline{x}}, \IK_{\underline{y}})} & {\Hom_{R\otimes R}(\Kyp(\IK_{\underline{x}}), \Kyp(\IK_{\underline{y}}))} \\
	{\Hom_{\Hecke}(\Delta_e,\IK_{\underline{z}}) } & {\Hom_{R\otimes R}(\Kyp(\Delta_e),\Kyp(\IK_{\underline{z}}))}
	\arrow["\Kyp", from=1-1, to=1-2]
	\arrow["\wr"', from=1-1, to=2-1]
	\arrow["\wr", from=1-2, to=2-2]
	\arrow["\Kyp", from=2-1, to=2-2]
\end{tikzcd}\]
where $\underline{z}$ is the concatenation of $\op{op}(\underline{x})$ and $\underline{y}.$ We hence need to show that the bottom vertical map in the commutative square is an isomorphism. 

By \Cref{lem:pureimpliespointwisepure}, the object $\IK_{\underline{z}}$ is pointwise $!$-pure. By \Cref{rem:pointwisepureyieldsfiltrations}, the object $\IK_{\underline{z}}$ admits a filtration whose associated graded are costandard objects, that is, there are objects $E^{k}\in \Hecke$ for $k=0,\dots,n$ with $E^0=0$ and $E^n=\IK_{\underline{z}}$ and cofiber sequences of the form
\[\begin{tikzcd}[column sep=scriptsize,row sep=tiny]
	{E^{k-1}} & {E^k} & {\nabla_{w_k}}
	\arrow[from=1-2, to=1-3]
	\arrow[from=1-1, to=1-2]
\end{tikzcd}\]
for some $w_k\in W.$ We show by induction on $k$ that the map 
\[\begin{tikzcd}
	{\Hom_{\Hecke}(\Delta_e,E^k) } & {\Hom_{R\otimes R}(R,\Kyp(E^k))}
	\arrow["\Kyp", from=1-1, to=1-2]
\end{tikzcd}\]
is an isomorphism for all $k$.
Consider the maps of exact sequences associated to the above cofiber sequence
\[\begin{tikzcd}
	& {\Hom(\Delta_e,E^{k-1})} & {\Hom(\Delta_e,E^{k})} & {\Hom(\Delta_e,\nabla_{w})} & 0 \\
	0 & {\Hom(R,\Kyp(E^{k-1}))} & {\Hom(R,\Kyp(E^{k}))} & {\Hom(R,R_{w})}
	\arrow[from=1-2, to=1-3]
	\arrow["\wr", from=1-2, to=2-2]
	\arrow[from=1-2, to=2-2]
	\arrow[from=1-3, to=1-4]
	\arrow[from=1-3, to=2-3]
	\arrow[from=1-4, to=1-5]
	\arrow["\wr", from=1-4, to=2-4]
	\arrow[from=1-4, to=2-4]
	\arrow[from=2-1, to=2-2]
	\arrow[from=2-2, to=2-3]
	\arrow[from=2-3, to=2-4]
\end{tikzcd}\]
where we abbreviate $w=w_k$ and the morphisms in the top row are in $\Hecke$ and those in the bottom row are in $\D(R\otimes R)$. We now explain why the diagram has the claimed form and hence by the five lemma the middle vertical arrow is an isomorphism.

The left vertical arrow is an isomorphism by induction.
For the right vertical arrow, we consider two cases.
If $w=e$, both source and target are isomorphic to $R$ and the map is clearly an isomorphism.
If $w\neq e$, the source of the arrow is zero by base change since $\Delta_e=i_{e,!}\un$ and $\nabla_w=i_{w,*}\un$ and $B\cap BwB=\emptyset$. Moreover, the target of the arrow is zero since $\Kyp(\Delta_e)=R$, $\Kyp(\nabla_w)=R_w$ and $\Hom_{R\otimes R}(R,R_w)=0$ by \Cref{lem:hombetweenstandardbimodules}. In this step we use that the Kac--Moody datum $\mathcal{D}$ is free.

The left-most term in the bottom row is $\Hom(R,R_{w}[-1])=0$, since both $R$ and $\Hom(R,R_{w})$ are in the heart of the standard $t$-structure on $\D(R\otimes R)$. The right-most term in the top row is $\Hom(\Delta_e,E^{k-1})=0$ using \Cref{lem:pointwisepureformality} and that $E^{k-1}$ is pointwise $!$-pure and $\Delta_e$ is pointwise $*$-pure.
\end{proof}
\begin{corollary} \label{cor:kyppuresoergelbimodules}
    Assume that the Kac--Moody datum $\Dd$ is free and of simply-connected type. Then the functor $\Kyp$ yields a monoidal equivalence between pure objects in the $K$-theoretic Hecke category and $K$-theory Soergel bimodules 
    $$\Heckepure\to \SBimK_\Dd$$
    mapping $E_s$ to $R\otimes_{R^s}R.$ 
\end{corollary}
\begin{remark}
We warn that the functor $\Kyp$ is not fully faithful on $\Hecke$, but just when restricted to pure objects. Already for $G=\Gm$, we have $\Hom_{\Hecke}(\Delta_e,\Delta_e[1])=0$ while $\Hom_{\D(R\otimes R)}(R,R[1])\neq 0$.
\end{remark}
\subsection{Soergel-theoretic description of the $K$-theoretic Hecke category} By combining the Erweiterungssatz, see \Cref{thm:erweiterungssatz} and \Cref{cor:kyppuresoergelbimodules}, and the formality result, see \Cref{cor:formalityofkmotivichecke}, we obtain the following `combinatorial' description of the $K$-theoretic Hecke category.
\begin{theorem}\label{thm:soergelktheoreticheckecategory}
    Assume that the Kac--Moody datum $\mathcal{D}$ is free and of simply-connected type, then there is an equivalence of monoidal categories
    $$\HeckeK_{\Dd}\to \Ch^b(\SBimK_{\Dd})$$
    between the $K$-theoretic Hecke category and the category of bounded chain complexes of $K$-theoretic Soergel bimodules associated to $\mathcal{D}.$ 
\end{theorem}
\section{Constructible monodromic sheaves}\label{sec:monodromicsheaves}
In this section we develop a theory of constructible monodromic sheaves.
\subsection{Sheaves on topological spaces}

Let $\Lambda$ be a regular Noetherian ring of finite global dimension. For a locally compact Hausdorff topological space $X$ we denote by $\D(X, \Lambda)$ the full derived category of sheaves of $\Lambda$-modules on $X$. We denote by $\Haus$ the category of locally compact Hausdorff topological spaces. 
\begin{theorem}[\protect{\cite[Lecture 7]{scholzeSixFunctorFormalisms}}]
There exists a $6$-functor formalism 
\begin{align*}
\D(-,\Lambda) : \Corr(\Haus) &\rightarrow \mathrm{Pr}_\Lambda \\
X &\mapsto \D(X,\Lambda).
\end{align*}
\end{theorem}

In particular, for all $X \in \Haus$ the category $\D(X,\Lambda)$ is a closed symmetric monoidal category. We denote the tensor product by $\otimes$ and the internal mapping spaces by $\iHom$. For all $f : X \to Y$, we have the usual functors $f^!, f^*, f_!$ and $f_*$ between $\D(X,\Lambda)$ and $\D(Y, \Lambda)$.

\begin{proposition}[\protect{\cite[Proposition 7.3]{scholzeSixFunctorFormalisms}}]
We equip the category $\Haus$ with the following Grothendieck topology. A collection of maps $(f_i : X_i \to X)_{i \in I}$ forms a cover if for any compact $K \subset X$ there exists a finite set $J \subset I$ and compacts $K_j \subset X_j$ for $j \in J$ such that $K = \cup_{j \in J} f_j(K_j)$. The functor $X \mapsto \D(X, \Lambda)$ is a sheaf in this topology. 
\end{proposition}

Using \cite[Proposition 4.17]{scholzeSixFunctorFormalisms} and \cite[Proposition A.5.16]{mannPAdic6FunctorFormalism2022} we can extend the $6$-functor formalism to stacks on the category $\Haus$. We give some example of such stacks, which we will use in the later sections. 
\begin{enumerate}
\item If $X = \varinjlim_i X_i$ is an increasing union of locally compact Hausdorff spaces where the transitions maps are closed immersions, then we have 
$$\D(X, \Lambda) = \varinjlim_{i} \D(X_i, \Lambda)$$
where the transitions maps are with respect to the $!$-pushforward along the inclusions. 
\item If $X = Y/G$ is a quotient then 
$$\D(X, \Lambda) = \varprojlim_{i} \D(Y^{\times_X i+1}, \Lambda),$$
where the transition maps are relative to the $*$-pullback and $Y^{\times_X i+1}=Y\times_X\dots\times_XY$ with $(i+1)$ copies of $Y$.
\end{enumerate}

\subsection{Three ways to present monodromic sheaves}

Let $H$ be a commutative Lie group with contractible universal cover which we denote by $\tilde{H}$. In particular $\tilde{H}$ sits in an extension 
\begin{equation*}
1 \rightarrow \pi_1(H, 1) \rightarrow \tilde{H} \rightarrow H \rightarrow 1. 
\end{equation*}

We denote ${\RH{H}} = \Lambda[\pi_1(H, 1)]$. There is a canonical map 
\begin{equation*}
\pi_1(H,1) \rightarrow (\RH{H})^{\times}, 
\end{equation*}
which defines a ${\RH{H}}$-local system $\Ll_H$ of rank one on $H$, called the free monodromic local system.

\begin{lemma}\label{lem:MultiplicativityFreeLocalSystem}
The sheaf $\Ll_H$ is canonically a multiplicative local system on $H$, that is there is a canonical isomorphism $m^*\Ll_H = \Ll_H \boxtimes_{{\RH{H}}} \Ll_H$ where $m : H \times H \rightarrow H$ is the multiplication map. 
\end{lemma} 

\begin{proof}
By definition of $\Ll_H$, if $f : H \rightarrow H'$ is a morphism of groups then we have a morphism ${\RH{H}} \rightarrow {\RH{H'}}$ and an isomorphism $f^*\Ll_{H'} = \Ll_H \otimes_{{\RH{H}}} {\RH{H'}}$. Applying this to the multiplication map which is a morphism since $H$ is commutative, we have 
\begin{align*}
m^*\Ll_H &= \Ll_{H \times H} \otimes_{\RH{H \times H}} {\RH{H}} \\
&= (\Ll_{H} \boxtimes_{\Lambda} \Ll_H) \otimes_{\RH{H} \otimes {\RH{H}}} {\RH{H}} \\ 
&= \Ll_H \boxtimes_{{\RH{H}}} \Ll_H. \qedhere
\end{align*}
\end{proof}

Given the data of a group $H$ together with a multiplicative local system $\Ll_H$, for all spaces $X$ with an action of $H$, there is a well-defined category of $(H, \Ll_H)$-equivariant sheaves on $X$. Let $X$ be a topological space with an action of $H$, we denote by $\D(X/(H,\Ll_H), {\RH{H}})$ the category of $(H,\Ll_H)$-equivariant sheaves on $X$ as defined in \cite[Section 1.5]{gaitsgoryLocalGlobalVersions2020}.  

Similarly, we denote by $\D(X/\tilde{H}, \Lambda)$ the category of sheaves on the stack $X/\tilde{H}$ which is the same as the category of $\tilde{H}$-equivariant sheaves on $X$. Since $\tilde{H}$ is contractible the forgetful functor $\D(X/\tilde{H}, \Lambda) \rightarrow \D(X, \Lambda)$ is fully faithful. 

\begin{theorem}\label{lemCharacterizationMonodromic}
Let $X$ be a topological space with an action of $H$. 
\begin{enumerate}
\item The forgetful functor $\D(X/(H, \Ll_H), {\RH{H}}) \rightarrow \D(X, \Lambda)$ forgetting both the equivariance and the ${\RH{H}}$-module structure along $\Lambda \rightarrow {\RH{H}}$ is fully faithful. 
\item The forgetful functor induces an equivalence 
\begin{equation*}
\D(X/\tilde{H}, \Lambda) = \D(X/(H,\Ll_H), {\RH{H}}).
\end{equation*}
\item Both categories are characterized as the full subcategories of sheaves on $X$ such that their $*$-pullback to $H$-orbits is locally constant. 
\item Let $A$ be a monodromic sheaf on $X$, then all $!$-pullbacks to $H$-orbits are locally constant. Conversely, if $A$ is constructible then $A$ is monodromic if and only if all of $!$-pullbacks to $H$-orbit are locally constant.
\end{enumerate}
\end{theorem} 

\begin{proof}
By definition of equivariant and twisted equivariant sheaves, their formation commutes with limits of categories and in particular satisfy descent. More precisely, if $X_* \rightarrow X$ is a simplicial $H$-resolution of $X$, then we have 
\begin{equation*}
\D(X/(H, \Ll_H), {\RH{H}}) = \varprojlim_{[n] \in \Delta^{op}} \D(X_n/(H, \Ll_H), {\RH{H}}).
\end{equation*}

Taking the bar resolution of $X$ provides us with a simplicial resolution of $X$ such that for all $[n] \in \Delta$ the action of $H$ on $X_n$ is free and the simplicial $H$-torsor $X_* \rightarrow X_*/H$ is trivial. Let $Y_*$ be a simplicial splitting, so that $X_* = Y_* \times H$. Then we have 
\begin{equation*}
\D(X_n/(H, \Ll_H), {\RH{H}}) = \D(Y_n, {\RH{H}}) \otimes_{\D({\RH{H}})} \D(H/(H, \Ll_H), {\RH{H}}).
\end{equation*}
Similarly, we have 
\begin{equation*}
\D(X/\tilde{H}, \Lambda) = \varprojlim_{\Delta^{op}} \D(X_n/\tilde{H}, \Lambda) = \varprojlim_{\Delta^{op}} \D(Y_n, \Lambda) \otimes_{\D(R)} \D(H/\tilde{H}, \Lambda).
\end{equation*}
Taking limits, we can assume first that the action of $H$ on $X$ is free and then that $X = H$. But now it is clear that there is an equivalence of categories $\D(H/(H, \Ll_H), {\RH{H}}) = \D({\RH{H}}) = \D(H/\tilde{H}, \Lambda)$. Hence, we get that the forgetful functor induces an equivalence $\D(X/(H, \Ll_H), {\RH{H}}) \rightarrow \D(X/\tilde{H}, \Lambda) \subset \D(X, \Lambda)$, this proves $(1)$ and $(2)$. 

Let us now show $(3)$. Consider the pair of adjoint functors $p_!: \D(X, \Lambda) \leftrightarrows \D(X/\tilde{H},\Lambda) : p^!$ where $p\colon X \to X/\tilde{H}$ is the quotient map. Since the pullback $p^!$ is conservative and continuous the category $\D(X/\tilde{H},\Lambda)$ is identified with the category of module over the monad $p^!p_!$. Since $p^!$ is fully faithful this monad is idempotent, that is the natural map $p^!p_!p^!p_! \to p^!p_!$ is an isomorphism and the category of modules over it is the full subcategory of objects $A \in \D(X, \Lambda)$ such that the map $A \to p^!p_!A$ is an isomorphism. So we have to show that $A$ satisfies this condition if and only if the $*$-pullback to all $H$-orbits are locally constant. 

For all $x \in X$, we denote by $o_x : H \to X, h \mapsto hx$, the orbit map of $x$. Let $A \in \D(X, \Lambda)$ be such that for all $x$, the object $o_x^*A$ is locally constant on $H$. Let $x \in X$ and consider the following Cartesian diagram 
\[\begin{tikzcd}
	H & {H/\tilde{H}} \\
	X & {X/\tilde{H}}
	\arrow["{p_x}", from=1-1, to=1-2]
	\arrow["{o_x}"', from=1-1, to=2-1]
	\arrow["{\tilde{o}_x}", from=1-2, to=2-2]
	\arrow["p"', from=2-1, to=2-2]
\end{tikzcd}\]
where the right column is the quotient by $\tilde{H}$ of the first column. We then have 
\begin{align*}
o_x^*A &= p_x^!p_{x,!}o_x^*A \\
&= p_x^!\tilde{o}_x^*p_{!}A \\
&=o_x^*p^!p_!A.
\end{align*}
Here, the first equality follows from the fact that $o_x^*A$ is locally constant on $H$. The second equality comes from base change and the last one from the smoothness of $\tilde{H}$. Hence, for all $x \in X$, the $*$-fiber at $x$ of the morphism $A \to p^!p_!A$ is an isomorphism, so the morphism itself is an isomorphism. Conversely, if the map $A \to p^!p_!A$ is an isomorphism, then pullback along $o_x$ shows via the same calculation that $o_x^*A$ is locally constant.

The case of $!$-fibers in $(4)$ is handled dually. Indeed, the functor $p^*$ is conservative and exhibits the category $\D(X/\tilde{H}, \Lambda)$ as the category of comodules over the comonad $p^*p_*$ (which is idempotent since $p^*$ is fully faithful). 
As before, if $A$ is monodromic, then the map $p^*p_*A \to A$ is an isomorphism. Applying the functor $o_x^!$ to this map shows that $o_x^!A$ is locally constant on $H$. 
Conversely, let $A$ be such that for all $x \in X$, the object $o_x^!A$ is locally constant and $A$ is constructible.  Then arguing as before, we get that the $!$-fiber at $x$ of the map $p^*p_*A \to A$ is an isomorphism and since $A$ is constructible, so is $p^*p_*A$ and the isomorphism can be detected on $!$-fibers. Hence, the map $p^*Ap_*A \to A$ is an isomorphism and $A$ is monodromic. 
\end{proof}

\begin{lemma}\label{lemmaFunctorialityMonodromic}
There is a well-defined four-functor formalism (so $f^*,f_*,f_!,f^!$) on stacks with $H$-action 
\begin{equation*}
X \mapsto \D(X/\tilde{H}, \Lambda) = \D(X/(H, \Ll_H), {\RH{H}}).
\end{equation*}
\end{lemma} 

\begin{proof}
It is simply a matter of noting that for all $f : X \rightarrow Y$ which is an $H$-equivariant map, the induced functors between $\D(X/\tilde{H}, \Lambda)$ and $\D(Y/\tilde{H}, \Lambda)$ are well-defined between subcategories of $\D(X, \Lambda)$ and $\D(Y,\Lambda)$ respectively. 
\end{proof}

Let $f : H_1 \to H_2$ be a surjective morphism of Lie groups which have contractible universal covers. Then we have a natural map $\RH{H_1} \to \RH{H_2}$. 

\begin{lemma}\label{lemmaHomogeneousSpace}
Let $X$ be a space with an action of $H_2$, there is an equivalence of categories 
\begin{equation*}
\D(X/\tilde{H}_1, \Lambda) = \D(X/\tilde{H}_2, \Lambda).
\end{equation*}
Moreover, we have 
\begin{equation*}
f_!(\Ll_{H_1} ) = \Ll_{H_2}[2 \dim(\ker(f))].
\end{equation*}
\end{lemma} 

\begin{proof}
For the first part of the lemma, we notice that both categories are full subcategories of $\D(X,R)$. It is therefore enough to check that they have the same objects. By \Cref{lemCharacterizationMonodromic}, they are the full subcategories of sheaves that are locally constant along $H_1$ (resp $H_2$) orbits. But since the map $H_1 \to H_2$ is surjective, the two conditions are equivalent. 

For the second part of the lemma, we first note that we have a commutative diagram 
\[\begin{tikzcd}
	{\ker \tilde{f}} & {\tilde{H}_1} & {\tilde{H}_2} \\
	{\ker f } & {H_1} & {H_2}
	\arrow[from=1-1, to=1-2]
	\arrow[from=1-1, to=2-1]
	\arrow["{\tilde{f}}", from=1-2, to=1-3]
	\arrow[from=1-2, to=2-2]
	\arrow[from=1-3, to=2-3]
	\arrow[from=2-1, to=2-2]
	\arrow["f"', from=2-2, to=2-3]
\end{tikzcd}\]
where $\tilde{f}$ is the map induced by $f$ on universal covers. Since both $\tilde{H}_1$ and $\tilde{H}_2$ are contractible it follows that $\ker(\tilde{f})$ is weakly contractible. Hence, we have $\tilde{f}_!\Lambda = \Lambda[-2\dim(\ker(\tilde{f}))] = \Lambda[-2\dim(\ker(f))]$, this implies the lemma.
\end{proof}

\subsection{Duality on monodromic sheaves}\label{subsectionDuality}

Let $X$ be a stack with an $H$-action. There are two natural Verdier type dualities that appear on free monodromic sheaves. 
\begin{remark}
On $\D(X/\tilde{H}, \Lambda)$ we can take $\Lambda$-linear Verdier duality. This one is not suited for applications since the free monodromic sheaf is not self dual for this duality as its stalks are not perfect complexes of $\Lambda$-modules. 
\end{remark}

On $\D(X/(H, \Ll_H), {\RH{H}})$, we can take ${\RH{H}}$-linear Verdier duality. However, the Verdier dual of a $(H, \Ll_H)$-equivariant sheaf is naturally $(H, \Ll_H^{\vee})$-equivariant (where $\Ll_H^{\vee}$ denotes the ${\RH{H}}$-dual of $\Ll_H$). Let $\iota : {\RH{H}} \rightarrow {\RH{H}}$ be the morphism induced by the inversion map of $H$. 

\begin{definition}
We define $\Dual' = {\RH{H}} \otimes_{\iota, {\RH{H}}} \Dual$. Since $\Ll_H \otimes_{{\RH{H}}, \iota} {\RH{H}} = \Ll_H^\vee$ the functor $\Dual'$ defines a duality functor of $\D(X/(H, \Ll_H), {\RH{H}})$. 
\end{definition}

\subsection{Constructible sheaves}

From now on we restrict our topological spaces to topological spaces of the form $X(\C)$ where $X$ is an algebraic variety of finite type over $\C$. A morphism $X(\C) \to Y(\C)$ is called algebraic if it comes from a morphism of algebraic varieties $X \to Y$. 

\begin{definition}\label{def:monodromicSheaves}
We denote by 
\begin{enumerate}
\item $\D_c(X, \Lambda)$ the category of constructible sheaves on $X$, this is the full subcategory of $\D(X,\Lambda)$ of sheaves $A$ for which there exists a Whitney stratification such that $A$ is locally constant with perfect $*$-stalk along all strata.
\item $\D_{\Ss-wc}(X, \Lambda)$ the category of weakly constructible sheaves on $X$, this the full subcategory of $\D(X,\Lambda)$ of sheaves $A$ for which there exists a stratification such that $A$ is locally constant. 
\end{enumerate}
For a quotient stack $X = Y/G$ the categories of (weakly-)constructible sheaves is the category of sheaves on $X$ whose pullback to $X$ is (weakly)-constructible. More generally for a category of twisted equivariant sheaves, the category of (weakly)-constructible twisted equivariant sheaves on $Y$ is the full subcategory of twisted equivariant sheaves such that the underlying sheaf on $Y$ is (weakly-)constructible. 
\end{definition}

\begin{remark}
Constructibility is defined in terms of $*$-stalks, the same definition using $!$-fibers is a priori not equivalent.
\end{remark}

\begin{theorem}[\cite{maximConstructibleSheafComplexes2022}]\label{thmMaxim}
Let $f : X \to Y$ be an algebraic morphism. Then the categories of constructible and weakly constructible sheaves are preserved under all functors $f_!,f_*, f^!, f^*$. 
\end{theorem}

\begin{definition}\label{def:monodromicconstructiblesheaves}
Let $H$ be a Lie group with contractible simply connected cover and $X$ be an $H$-space. Then we define the category of \emph{constructible monodromic sheaves} to be 
\begin{equation*}
\D_c(X, \Lambda)_{\mon} = \D_{\Ss-wc}(X/\tilde{H}, \Lambda) \cap \D_c(X/(H, \Ll), {\RH{H}})
\end{equation*}
of sheaves that are $\Lambda$-weakly constructible but have perfect ${\RH{H}}$-stalks. 
\end{definition}

\begin{lemma}
For all $H$-equivariant algebraic morphisms $f : X \to Y$ the category of constructible monodromic sheaves is preserved under all functors $f_!,f_*,f^*,f^!$.  
\end{lemma} 

\begin{proof}
This is an immediate combination of \Cref{thmMaxim} and \Cref{lemmaFunctorialityMonodromic}.
\end{proof}

The next proposition is an immediate consequence of standard theorems for Verdier duality. 
\begin{proposition}
For $f : X \to Y$ an $H$-equivariant map, we have $\Dual'f^! = f^*\Dual'$ and $\Dual'f_! = f_*\Dual'$. Moreover, for an ${\RH{H}}$-constructible sheaf $A$, the biduality map $A \to \Dual'\Dual'(A)$ is an isomorphism. 
\end{proposition}

\section{Soergel description of monodromic sheaves on extended flag varieties}\label{sec:soergelmonodromic}
In this section we will give a Soergel-theoretic description of the universal monodromic Hecke category $\Heckemon_{\Dd}$ which consists of $T\times T$-monodromic sheaves on the stack $U\backslash G/U$ of a Kac--Moody group $G$ associated to a Kac--Moody root datum $\Dd$. In this section, we assume that the root datum $\Dd$ is cofree and of adjoint type as explained in \Cref{sec:kacmoodyrootdatum}. We also fix a set of lifts $(\dot{w})$ of elements of $W$ in $G$. We abbreviate the group algebra of the fundamental group of the torus by $\LDR=\Z[\pi_1(T)].$

\subsection{Monodromic sheaves on flag varieties}

Denote by $G\supset B\supset T$ the Kac--Moody group together with a Borel subgroup and maximal torus associated to the Kac--Moody root datum $\Dd.$ Denote by $U\subset B$ the unipotent radical, by $W\supset S$ the Weyl group and simple reflections. In this section, we identify all the groups with their complex points equipped with the analytic topology. 

On the stack $U \backslash G/U$, there are three actions of tori that we can consider.
\begin{enumerate}
\item The action of $T$ acting by left translations, 
\item the action of $T$ acting by right translations and
\item the action of $T \times T$ acting by left and right translations. 
\end{enumerate}
We denote by $\Hecke^{\lef}, \Hecke^{\rig}$ and $\Hecke^{\lef\rig}$ the corresponding categories of constructible monodromic sheaves $\D_c(U\bs G/U)_\mon$, see \Cref{def:monodromicconstructiblesheaves}. There are obvious forgetful functors $\Hecke^{\lef} \xleftarrow{\oblv^{\rig}} \Hecke^{\lef\rig} \xrightarrow{\oblv^{\lef}} \Hecke^{\rig}$.

\begin{lemma}
Both functors $\oblv^{\rig}$ and $\oblv^{\lef}$ are equivalences. 
\end{lemma}

\begin{proof}
Firstly, using the Bruhat decomposition of $U\backslash G/U$, we see that the categories $\Hecke^{\lef}, \Hecke^{\rig}$ and $\Hecke^{\lef\rig}$ are obtained by gluing the categories of monodromic sheaves on each stratum. Since the forgetful functors are compatible with the four functors, they are compatible with the gluing, in particular, we only need to check the statement on each stratum, this is a direct application of \Cref{lemmaHomogeneousSpace}. 
\end{proof}

\begin{definition}
The \emph{universal monodromic Hecke category} $\Hecke=\Heckemon_{\Dd}$ is defined as either one of the three categories $\Hecke^{\lef}, \Hecke^{\rig}$ or $\Hecke^{\lef\rig}$. 
\end{definition}
We now introduce standard, costandard and tilting perverse sheaves in $\Hecke$.
Let $n \in N(T)$ be an element in the normalizer of $T$ and denote by $w$ its image in the Weyl group of $G$. The choice of $n$ yields a decomposition of the Bruhat stratum $BwB \simeq U \times T \times U_w$ given by $(u,t,v) \mapsto utnv$ where $U_w = U \cap \Ad(n)U$. We denote by $\nu_n : BwB \rightarrow T$ the projection onto $T$. This map is equivariant for the action of $U \times U$ acting trivially on the target, hence we have a well-defined map of stacks $\nu_n : U \backslash BwB/U \to T$ which we will also call $\nu_n$. 
\begin{definition}
Let $n \in N(T)$, we denote by 
\begin{enumerate}
\item $\Delta_n = i_{w,!} \nu_n^*\Ll_T[\ell(w) + \dim(T)]$ the standard sheaf, 
\item $\nabla_n = i_{w,*} \nu_n^*\Ll_T[\ell(w) + \dim(T)]$ the costandard sheaf. 
\end{enumerate}

\end{definition}
The standard and costandard sheaves are monodromic constructible sheaves (using either the left or right $T$-action) as the sheaf $\Ll_T$ is an $\LDR$-constructible sheaf on $T$. Since the inclusions $i_w$ are affine, they are perverse. 

\begin{remark}
Up to isomorphism, the sheaves $\Delta_n$ depend only on $w$. In what follows, we will denote by $\Delta_w = \Delta_{\dot{w}}$ where $\dot{w}$ is the lift of the element $w$ we have chosen. 
\end{remark}

\begin{definition}
Let $A \in \Hecke$ be a perverse sheaf,
\begin{enumerate}
\item a \emph{$\Delta$-flag} for $A$ is a filtration such that the graded pieces are isomorphic to standard sheaves,
\item a \emph{$\nabla$-flag} for $A$ is a filtration such that the graded pieces are isomorphic to costandard sheaves and
\item the sheaf $A$ is \emph{tilting} if it has both a $\Delta$-flag and a $\nabla$-flag. 
\end{enumerate}
We denote by $\Hecketilt\subset \Hecke$ the full subcategory of tilting sheaves.
\end{definition}

\subsection{Monoidal Structure}

We equip the category $\Hecke$ with a convolution structure. Consider the diagram 
\[\begin{tikzcd}
	& {U \backslash G \times^U G/U} & {U \backslash G/U} \\
	{U \backslash G/U} && {U \backslash G/U}
	\arrow["m", from=1-2, to=1-3]
	\arrow["{p_2}"', from=1-2, to=2-3]
	\arrow["{p_1}", from=1-2, to=2-1]
\end{tikzcd}\]
The convolution product is defined in two steps. Let $A, B \in \Hecke$ which we see as the category of $(T \times T, \Ll_{T \times T})$-equivariant sheaves on $U \backslash G/U$, then $p_1^*A \boxtimes_{\mathbb{Z}} p_2^*B$ is naturally a sheaf of $\LDR^{\otimes 4}$-modules, and it is constructible as an $\LDR^{\otimes 4}$-sheaf on $U \backslash G \times^U G/U$. We define the convolution of $A$ and $B$ as 
\begin{equation*}
A \conv B = \For_{\LDR^{\otimes 2}}^{\LDR^{\otimes 4}}m_!(p_1^*A \otimes_{\mathbb{Z}} p_2^*B)[\dim T] 
\end{equation*}
where the functor $\For_{\LDR^{\otimes 2}}^{\LDR^{\otimes 4}}$ is the forgetful functor induced by the inclusion $\LDR \otimes \LDR \rightarrow \LDR^{\otimes 4}$ induced by the outer inclusions. Keeping track of the actions of $T$, there is a canonical $\LDR \otimes \LDR$-linear isomorphism 
\begin{equation*}
A \conv B = m_!(p_1^*\For^{\rig}A \otimes_{\mathbb{Z}} p_2^*\For^{\lef}B)[\dim T].
\end{equation*}
In particular, since both $\For^{\rig}A$ and $\For^{\lef}B$ are $\LDR$-constructible, the sheaf $A \conv B$ is $\LDR \otimes \LDR$-constructible. 
There is also a dual convolution defined as
\begin{equation*}
	A \conv^! B = m_*(p_1^!\For^{\rig}A \otimes_{\mathbb{Z}} p_2^!\For^{\lef}B).
\end{equation*}
\subsection{Duality}
Recall that we introduced the functor $\Dual'$ as a replacement of Verdier duality for monodromic sheaves, see \Cref{subsectionDuality}. For $A \in \Hecke$, we introduce the following duality functor 
\begin{equation*}
	\Dual^{-}(A) = \mathrm{inv}^*\Dual'(A)(\varepsilon)[-\dim T]
\end{equation*}
where $(\varepsilon)$ is again the twist of the $\LDR$-structure by induced by the inversion map. 
Let us collect the following standard properties of the dualities.
\begin{lemma}\label{lem:propertiesdualitymonodromic}
	On monodromic Hecke category, we have the following statements.
		\begin{enumerate}
			\item $(\Dual')^2\cong \id\cong (\DualI)^2$ on $\Hecke.$
		    \item There are isomorphisms $\Dual'(\Delta_n)\cong \nabla_n$ and $\DualI(\Delta_n)\cong \nabla_{n^{-1}}.$
			\item There is a natural isomorphism $\Dual(A)\conv^!\Dual(B)\to \Dual(A\conv B)$ and $\Dual(A\conv^! B)\to \Dual(A)\conv\Dual(B).$
			\item There is a natural isomorphism $\Dual^-(A\conv B)\to \Dual^-(B)\conv^!\Dual^-(A)$.
		\end{enumerate}
\end{lemma}
\begin{lemma}
For $A, B, C \in \Hecke$, there are canonical isomorphisms
$$\Hom(A \conv B, C) = \Hom(A, C \conv^! \Dual^{-}(B)) = \Hom(B, \Dual^{-}(A) \conv^! C).$$
\end{lemma}
\begin{proof}
This can be shown by a yoga of adjunctions as in \Cref{lem:dualizability}.
\end{proof}
\begin{lemma}\label{lemConvolutionMonodromic}
	Let $n,n' \in N(T)$ and assume that $\ell(n) + \ell(n') = \ell(nn')$, then 
	\begin{enumerate}
	\item $\Delta_n \conv \Delta_{n'} = \Delta_{nn'}$,
	\item $\nabla_n \conv^! \nabla_{n'} = \nabla_{nn'}$,
	\end{enumerate} 
	\end{lemma} 
	\begin{proof}
		Follows from the fact that $BnB \times^B Bn'B \to BnnB$ is an isomorphism.
\end{proof}
\begin{theorem} \label{thm:rigiditymonodromic}
There is a natural equivalence $\conv\simeq \conv^!$ and all objects in $\Hecke$ are left and right dualisable with left and right duals canonically identified with $\Dual^{-}(-)$. 
\end{theorem}
\begin{proof}
	By \Cref{lemInversion}, for a simple reflection $s$ the object $\Delta_s$ is dualisable with respect to $\conv$. Hence, so is $\Delta_n$ for all $n\in N(T)$ using \Cref{lemConvolutionMonodromic}. Now we can argue as in \Cref{thm:rigiditykmotivichecke}, since the objects $\Delta_n$ generate the category.
\end{proof}

\begin{lemma}
	We have $\Delta_n \conv \nabla_{n^{-1}} = \nabla_{n^{-1}} \conv \Delta_n = \Delta_1 = \nabla_1$.
\end{lemma}
\begin{proof}
	This follows from \Cref{lemConvolutionMonodromic}, \Cref{thm:rigiditymonodromic} and the case of simple reflection which is \Cref{lemInversion}.
\end{proof}

\begin{remark}
Note that any standard or costandard sheaf is invertible, in particular we get that $\Dual^{-}(\Delta_n) = \nabla_{n^{-1}}$ and $\Dual^{-}(\nabla_{n}) = \Delta_{n^{-1}}$ for all $n \in N(T)$ and that $\Dual^{-}$ preserve the category of tilting sheaves. 
\end{remark}

\subsection{Rank one calculations}\label{sec:rankonecalculations}

The goal of this section is to do the rank one calculations. Most proofs will reduce to arguments of \cite{taylorUniversalMonodromicTilting2023}. In \emph{loc.~cit.}, the author assumes the group $G$ to be adjoint. From now on, we assume that $G$ is of semisimple rank one of adjoint type so it has a connected center. 

Recall that $\LDR = \Z[X_*(T)]$ is the group ring on the space of cocharacters of $T$. For $\lambda \in X_*$, we denote by $e^{\lambda} \in \LDR$ the corresponding element. We denote by $\alpha^{\vee}$ the simple coroot of $G$ and by $\alpha$ the simple root of $G$. We also fix as in \Cref{subsectionSoergelBimodules} a cocharacter $\varpi^{\vee}$ such that $\langle \varpi^{\vee}, \alpha \rangle = 1$.

We now consider the category $\Hecke$ for $G$ and with the objects $\Delta_1 = \nabla_1$, $\Delta_s$ and $\nabla_s$ (note that the last two objects are defined up to isomorphism).  
\begin{lemma}\label{lemInversion}
There is an isomorphism 
$$\Delta_s \conv \nabla_s = \nabla_s \conv \Delta_s = \Delta_1.$$
\end{lemma}
\begin{proof}
The argument of \cite[Proposition A.4]{taylorUniversalMonodromicTilting2023} holds in our context.
\end{proof}

\begin{lemma}\label{lemStructureTS}
There exists a tilting object $T_s$ satisfying the following 
\begin{enumerate}
\item There are two short exact sequences 
$$0 \to \Delta_s \to T_s \to \Delta_1 \to 0$$ 
and 
$$0 \to \nabla_1 \to T_s \to \nabla_s \to 0.$$
\item Any object $K$ satisfying $(i)$ is isomorphic to $T_s$.
\item There is an isomorphism $\mathrm{End}(T_s) = \LDR \otimes_{\LDR^s} \LDR$. 
\end{enumerate}
\end{lemma}

\begin{proof}
(1) The first point is shown as in \cite[Lemma 5.1]{taylorUniversalMonodromicTilting2023}, namely we have (noncanonical) isomorphisms
$$\Ext^1(\Delta_1, \Delta_s) = \Ext^1(\nabla_s, \nabla_1) = \LDR/(1 - e^{\alpha^{\vee}}).$$
The element $1$ then yields an object $T_s$ with the two desired short exact sequences. 

(2) We can argue as in \cite[Lemma 5.2]{taylorUniversalMonodromicTilting2023}. The extension $K$ corresponds to an element $a \in \LDR/(1 - e^{\alpha^{\vee}})$, it is then enough to check that this element is invertible in this ring. As the ring $\LDR/(1 - e^{\alpha^{\vee}})$ is reduced by \Cref{lemRingReduced}, it is enough to show that the image of $a$ in all residue fields of $\LDR/(1 - e^{\alpha^{\vee}})$ are nonzero. 

Let $\mathfrak{p}$ be a prime ideal of $\LDR/(1 - e^{\alpha^{\vee}})$, which we consider as a prime ideal of $\LDR$. Tensoring the first exact sequence with $\otimes_{\LDR} \LDR/\mathfrak{p}$ remains exacts as $\Delta_1 \otimes_{\LDR} \LDR/\mathfrak{p}$ lies in perverse degree $0$. Hence, we have a short exact sequence 
$$0 \to \Delta_s \otimes_{\LDR} \LDR/\mathfrak{p} \to K \otimes_{\LDR} \LDR/\mathfrak{p} \to \Delta_1 \otimes_{\LDR} \LDR/\mathfrak{p} \to 0$$
which corresponds to the image of $a$ along the induced map $\Ext^1(\Delta_s, \Delta_1) \to \Ext^1(\Delta_s \otimes_{\LDR} \LDR/\mathfrak{p}, \Delta_1 \otimes_{\LDR} \LDR/\mathfrak{p}) = \LDR/(1-e^{\alpha^{\vee}}) \otimes_{\LDR} \LDR/\mathfrak{p} = \LDR/\mathfrak{p}$.

If we suppose that this extension is split, that is $a$ is $0$ in $\LDR/\mathfrak{p}$, then we can find a splitting $K\otimes_{\LDR} \LDR/\mathfrak{p} = \Delta_1 \otimes_{\LDR} \LDR/\mathfrak{p} \oplus \Delta_s \otimes_{\LDR} \LDR/\mathfrak{p}$. Using the second short exact sequence, we get
$$0 \to \nabla_1 \otimes_{\LDR} \LDR/\mathfrak{p} \to  \Delta_1 \otimes_{\LDR} \LDR/\mathfrak{p} \oplus \Delta_s \otimes_{\LDR} \LDR/\mathfrak{p} \to \nabla_s \otimes_{\LDR} \LDR/\mathfrak{p} \to 0.$$
As $\Hom(\Delta_1 \otimes_{\LDR} \LDR/\mathfrak{p}, \nabla_s \otimes_{\LDR} \LDR/\mathfrak{p}) = 0$, the second map in the short exact sequence induces an isomorphism $\Delta_s \otimes_{\LDR} \LDR/\mathfrak{p} = \nabla_s \otimes_{\LDR} \LDR/\mathfrak{p}$ however this contradicts the assumption that $\mathfrak{p}$ was an ideal of $\LDR/(1 - e^{\alpha^{\vee}})$, hence $a$ is nonzero. 

(3) We argue as in \cite[Proposition 5.4]{taylorUniversalMonodromicTilting2023}. Firstly the monodromy map and the functors $\mathrm{gr}$ of \cite[6.3]{bezrukavnikovTopologicalApproachSoergel2020} yield maps
$$\LDR \otimes_{\LDR^s} \LDR \to \mathrm{End}(T_s) \to \LDR_e \oplus \LDR_s.$$
We want to show that the first map is an isomorphism. Their composition is injective by \Cref{lemInjectivity} hence the first map is injective. The injectivity of the second map follows from \cite[Corollary 6.3]{bezrukavnikovTopologicalApproachSoergel2020}. It is then enough to show that both algebras have the same image in $\LDR_e \oplus \LDR_s$. Let $a \in \mathrm{End}(T_s)$ and denote by $a_e$ and $a_s$ its components in $\LDR_e \oplus \LDR_s$. Since $\LDR \otimes_{\LDR^s} \LDR \to \LDR_e \oplus \LDR_s \to \LDR_e$ is surjective by \Cref{lemInjectivity}, we can find $b \in \LDR \otimes_{\LDR^s} \LDR$ such that $b_e = -a_e$. We consider the endomorphism $a + b$ of $T_s$. By the same Lemma, it is enough to show that $a_s + b_s$ lie in the ideal $(1 - e^{\alpha})$. Equivalently, it is enough to show that this element is zero in $\LDR/(1 - e^{\alpha})$ and since this ring is reduced, it is enough to show that this element vanishes in all residue fields of this ring. But the condition on $\varpi$ implies that after reducing modulo some prime ideal containing $1 - e^{\alpha}$, the extension 
$$0 \to \Delta_s \to T_s \to \Delta_1 \to 0$$
is not split. As $a_e + b_e = 0$, the map $a + b : T_s \to T_s$ factors through $\Delta_s$ and after reducing modulo some ideal it can only be zero only if the extension splits.
\end{proof}

\begin{remark}\label{remarkConvolTs}
As observed in \cite[Lemma 5.3]{taylorUniversalMonodromicTilting2023}, it follows from \Cref{lemStructureTS} by convolving the second exact sequence with $\Delta_s$ that there is a triangle 
$$\Delta_s \to \Delta_s \conv T_s \to \Delta_s \conv \nabla_s.$$
Since $\Delta_s \conv \nabla_s = \Delta_1$, this triangle is then a short exact sequence of perverse sheaves and exhibits that we have an isomorphism $$\Delta_s \conv T_s = T_s.$$
\end{remark}

\subsection{Tilting sheaves}

We now collect some results about convolution, duality and the existence of perverse tilting sheaves.

For this, denote by $\Hecke^s\subset \Hecke$ the Hecke category for the minimal parabolic corresponding to a simple reflection $s\in S.$ Then, by the discussion in \Cref{lemStructureTS}, there is an indecomposable tilting sheaf $T_s\in \Hecke^s$ such that the multiplicity of $\Delta_s$ in a $\Delta$-flag of $T_s$ is $1$. The sheaf $T_s$ is unique up to isomorphism and in this section we fix an arbitrary representative of this sheaf in its isomorphism class. Once we introduce the $\mathbb{V}$-functor, we will have a canonical way of choosing such a representative.

\begin{lemma}\label{lemmaConvolutionTilting}
Let $T,T' \in \Hecketilt$ be two tilting perverse sheaves, then $T \conv T'$ is tilting and perverse. 
\end{lemma} 

\begin{proof}
We mimic the argument of \cite[Proposition 4.3.3]{bezrukavnikovKoszulDualityKacMoody2013}. We first show that for $w,w' \in W$ we have 
\begin{align*}
\Delta_w \conv \Delta_{w'} &\in \langle \Delta_v[\leq 0], v \in W \rangle \subset \Hecke \text{ and}\\
\nabla_w \conv \nabla_{w'} &\in \langle \nabla_v[\geq 0], v \in W \rangle \subset \Hecke,
\end{align*}
where $\langle - \rangle$ denotes the full subcategory generated by extensions. 
Let us prove the first point, the second one is proven in a similar way. We can argue by induction on $\ell(w)$ and reduce to the case $w = s$ is a simple reflection. Then if $\ell(sw') = \ell(w') + 1$, we have $\Delta_s \conv \Delta_{w'} = \Delta_{sw'}$. If $\ell(sw') = \ell(w') - 1$, then we use the short exact sequence 
\begin{equation*}
0 \to \Delta_s \to T_s \to \Delta_1 \to 0, 
\end{equation*}
from \Cref{lemStructureTS}. After convolving with $\Delta_s$, we get a triangle
\begin{equation*}
\Delta_s \conv \Delta_s \to \Delta_s \conv T_s \to \Delta_s 
\end{equation*}
hence $\Delta_s \conv \Delta_s  \in \langle \Delta_s \conv T_s, \Delta_s[-1] \rangle$. By \Cref{remarkConvolTs}, we have $\Delta_s \conv T_s \simeq T_s$ hence $\Delta_s \conv \Delta_{w'}$ is a successive extension of $\Delta_{w'}[-1], \Delta_{sw'}$ and $\Delta_{w'}$.

The sheaf $T \conv T'$ now belongs to the category $\langle \Delta_v[\leq 0], v \in W \rangle \cap \langle \nabla_v[\geq 0], v \in W \rangle$. The rest of the argument of \cite{bezrukavnikovKoszulDualityKacMoody2013} follows verbatim. 
\end{proof}
\begin{lemma}
	There is an isomorphism 
    $\Dual^{-}(T_s) \cong T_s.$
\end{lemma} 
\begin{proof}
	The two triangles defining $T_s$, see \Cref{lemStructureTS}, are exchanged by $\Dual^-$, using that $\Dual^-(\Delta_s)=\nabla_s$ and $\Dual^-(\Delta_e)=\nabla_e.$

	By \Cref{lemStructureTS} we get $\Dual^-(T_s) = T_s$.
\end{proof}	

\begin{remark}
Note that we have used here that the root datum is of adjoint type. 
\end{remark}

\begin{lemma}\label{lemmaExistenceTilting}
For all $w$, there exists a tilting sheaf supported on the closure of $BwB$ and such that the multiplicity of $\Delta_w$ in any of its $\Delta$-flags is one. 
\end{lemma}

\begin{proof}
	Write $w = s_1\dots s_n$ be a reduced expression for $w$, then by \Cref{lemmaConvolutionTilting} the sheaf $T = T_{s_1} \conv \dots \conv T_{s_n}$ is tilting and supported on the closure of $BwB$. Since the map $Bs_1B \times^B Bs_2B \times^B \dots Bs_nB \to BwB$ is an isomorphism the multiplicity of $\Delta_w$ in a $\Delta$-flag of $T$ is $1$. 
\end{proof}

\begin{corollary}
	The category of perverse tilting sheaves is generated by the objects $T_s$ as a monoidal, additive and idempotent closed category
$$\Hecketilt=\genbuild{T_s}{s\in S}_{\conv, \inplus,\cong}\subset \Hecke.$$
\end{corollary}

\begin{corollary}\label{cor:tiltinggeneratehecke}
	The Hecke category is stably generated by the tilting perverse sheaves
$$\genbuildexplicit{\Hecketilt}_{\op{stb}}=\Hecke.$$
\end{corollary}

\subsection{Formality} We now prove the formality of the universal monodromic Hecke category using tilting objects.
We first recall the following standard property of tilting objects in our context.
\begin{lemma}
	Let $T,T'\in \Hecketilt$. Then $\Hom_\Hecke(T,T'[n])=0$ if $n\neq 0$
\end{lemma}
\begin{proof}
	Using that $T$ has a $\Delta$-flag and $T'$ has a $\nabla$-flag, the statement follows by induction and $\Hom_\Hecke(\Delta_w,\nabla_w'[n])=0$ for all $w,w'\in W$ and $n\neq 0.$
\end{proof}
Using that the tilting objects generate the Hecke category, see \Cref{cor:tiltinggeneratehecke}, \Cref{prop:weightcomplexequivalence} implies the following formality statement.
\begin{corollary}\label{corGenerationBottSamelsonMonodromic}
There is an equivalence of monoidal categories between the universal monodromic Hecke category and the category of bounded chain complexes of perverse tilting sheaves
    $$\Hh\stackrel{\sim}{\to}\Ch^b(\Hecketilt).$$
\end{corollary} 
\begin{remark}
	The equivalence in \Cref{corGenerationBottSamelsonMonodromic} uses the weight complex functor for the `tilting weight structure'. The inverse of the equivalence is given by
	$$\Ch^b(\Hecketilt)\stackrel{\sim}{\to}\D^b(\Hecke^{\heartsuit})\stackrel{\sim}{\to}\Hecke$$
	where the second equivalence is Beilinson's realization functor. We also refer to \cite[1.5 Proposition]{beilinsonTiltingExercises2004a} where a similar statement for constructible sheaves on $G/B$ is discussed.
\end{remark}

\subsection{The $\mathbb{V}$-functor}

We now define the $\V$-functor. Consider the following character of $U^{-}$
\begin{equation*}
\chi : U^{-} \rightarrow \prod_{\alpha \in \Delta} \Ga \xrightarrow{\Sigma} \Ga. 
\end{equation*}
Denote by $1 : \point \rightarrow U^{-}$ the inclusion of the point $1$ and $i : U^{-} \rightarrow U \backslash G/U$ the map induced by the inclusion of $U^{-}$. The $\mathbb{V}$-functor is defined as 
\begin{align*}
\Hecke &\rightarrow \D(\Z) \\
A &\mapsto 1^*\phi_{\chi}(i^*A),
\end{align*}
where $\phi_{\chi}$ denotes the vanishing cycle functor $\D(U^{-}) \rightarrow \D(\chi^{-1}(0))$. Using the presentation of the category $\Hecke$ as a category of equivariant $\LDR \otimes \LDR$-sheaves, the $\V$-functor factors through the category $\D(\LDR \otimes \LDR)$. Moreover, since the vanishing cycle functor is $t$-exact, it follows that the $\V$-functor is $t$-exact. 

Following \cite[Section 2.2.4]{liFunctionsCommutingStack2023}, there is a canonical lax-monoidal structure on $\V$ constructed as follows. Consider the map 
$$i^{(2)} : U^- \times U^- \to U \backslash G \times^U G/U$$
induced by the product of the two copies of the map $U^- \to G$. Consider the character 
$$\chi + \chi : U^- \times U^- \to \A^1,$$ 
and define 
$$\V_{U \backslash G \times^U G/U}(-) = 1^*\phi_{\chi+\chi}i^{(2),*}$$
where $\phi_{\chi + \chi}$ is the vanishing cycle functor with respect to $\chi + \chi$. Let us summarize all the objects in the following diagram 
\[\begin{tikzcd}
	{U \backslash G/U \times U\backslash G/U} & {U^- \times U^-} & {\A^1} \\
	{U \backslash G \times^{U} G/U} & {U^- \times U^-} & {\A^1} \\
	{U \backslash G/U} & {U^-} & {\A^1}
	\arrow[from=1-2, to=1-1]
	\arrow["{\chi + \chi}", from=1-2, to=1-3]
	\arrow["p", from=2-1, to=1-1]
	\arrow["m"', from=2-1, to=3-1]
	\arrow[equals, from=2-2, to=1-2]
	\arrow[from=2-2, to=2-1]
	\arrow["{\chi + \chi}", from=2-2, to=2-3]
	\arrow["m"', from=2-2, to=3-2]
	\arrow[equals, from=2-3, to=1-3]
	\arrow[equals, from=2-3, to=3-3]
	\arrow[from=3-2, to=3-1]
	\arrow["\chi"{description}, from=3-2, to=3-3]
\end{tikzcd}\]

Consider now the composition 
\begin{align*}
\V(-) \otimes_{\Z} \V(-) &= 1^*\phi_{\chi}(i^*(-)) \otimes 1^*\phi_{\chi}(i^*(-)) \\
&= 1^*(\phi_{\chi}(-) \boxtimes \phi_{\chi}(-)) \\
&\xrightarrow{\sim} 1^*(\phi_{\chi + \chi}((i^*-) \boxtimes (i^*(-)) \\
&= 1^*\phi_{\chi + \chi}(i^{(2),*}p^*(- \boxtimes -)) \\
&= \V_{U \backslash G \times^U G/U}(p^*( - \boxtimes -)). 
\end{align*}
In the above composition, the only nontrivial map is the Thom-Sebastiani map 
$$1^*\phi_{\chi}(-) \otimes 1^*\phi_{\chi}(-) \to 1^*\phi_{\chi + \chi}(- \boxtimes -),$$
which is a Künneth map for vanishing cycles and proven to be an isomorphism in full generality in \cite{ThomSebastianiThoeremMassey}. 

We now get, for all $A, B \in \Hecke$
\begin{align*}
\V(A) \otimes_{\Z} \V(B) &\simeq \V_{U \backslash G \times^U G/U}(p_1^*(A) \otimes_{\Z} p_2^*(B)) \\
&\rightarrow \V_{U \backslash G \times^U G/U}(m^!m_!(p_1^*(A) \otimes_{\Z} p_2^*(B)) = \V(A \conv B)
\end{align*}
where the last isomorphism comes from the compatibility of vanishing cycles with smooth pullback along the map $m$.

\begin{theorem}[\protect{\cite[Proposition 2.1]{taylorUniversalMonodromicTilting2023}}]
The functor $\V$ is equipped with a canonical monoidal structure $\Hecke \to \D(\LDR \otimes \LDR)$ where the target category is equipped with the convolution of bimodule structure, in particular, for $A, B \in \Hecke$, the map 
$$\V(A) \otimes_{\Z} \V(B) \to \V(A \conv B)$$
factors canonically through $\V(A) \otimes_{\LDR} \V(B).$
\end{theorem}

\begin{remark}
In \emph{loc.~cit.}, the proof is done in the setting of finite dimensional groups, but the same proof extends verbatim to the Kac--Moody setting. 
\end{remark}

Let $s$ be a simple reflection. By \cite{taylorUniversalMonodromicTilting2023}, the functor $\V$ restricted to $\Hh^s$ is represented by the object $T_s$. This property pins the object $T_s$ uniquely. 
\subsection{Struktursatz} 
We now prove the analogue of Soergel's Struktursatz \cite{soergelKategorieMathcalPerverse1990} in the setting of the universal monodromic Hecke category.
\begin{theorem}[Struktursatz]\label{thmStrukturSatzMonodromic}
	Recall that $\Dd$ is cofree and of adjoint type. Then the functor $\V: \Hecketilt\to \D(\LDR\otimes \LDR)$ is fully faithful.
	\end{theorem}
	
	\begin{proof}
	Using the discussion in \Cref{sec:rankonecalculations}, we have $\V(\nabla_s)=\Hom(T_s,\nabla_s)=R_s$ so by the monoidality of $\V$ we obtain that $\V(\nabla_w)\cong R_w$. 
	Now, let $T,T'\in \Hecketilt$. By the monoidality of $\V$ we get the commutative diagram
	
\[\begin{tikzcd}
	{\Hom(T',T)} & {\Hom_{\LDR\otimes \LDR}(\V(T'),\V(T))} \\
	{\Hom(\Delta_e,\Dual^-(T')\star T)} & {\Hom_{\LDR\otimes \LDR}(R,\V(\Dual^-(T')\star T)).}
	\arrow[from=1-1, to=1-2]
	\arrow["\wr", from=1-1, to=2-1]
	\arrow["\wr"', from=1-2, to=2-2]
	\arrow[from=2-1, to=2-2]
\end{tikzcd}\]
Hence, we are reduced to the case that $T'=\Delta_e$. We show the stronger statement that $\V$ yields an isomorphism
\begin{equation*}
	\Hom(\Delta_1, T) \to \Hom_{\LDR \otimes \LDR}(R, \V(T)). 
	\end{equation*}
	for objects $T$ with a $\nabla$-flag.
	We prove this by induction on the filtration length of $T$. If $T=\nabla_w$, then this is immediate since $\V(\nabla_w)=R_w$ and $\Hom_{\LDR\otimes\LDR}(R,R_w)=0$ if $w\neq 1$ using \Cref{lem:hombetweenstandardbimodules} and the assumption that $\mathcal{D}$ is cofree.
	For the induction step, pick a short exact sequence
	$T'\to T\to \nabla_w$. Using that there are no extension between standard and costandard objects we obtain the following diagram of exact sequences
	
\[\begin{tikzcd}
	0 & {\Hom(\Delta_1,T')} & {\Hom(\Delta_1,T)} & {\Hom(\Delta_1,\nabla_w)} & 0 \\
	0 & {\Hom(R,\V(T'))} & {\Hom(R,\V(T))} & {\Hom(R,R_w)}
	\arrow[from=1-1, to=1-2]
	\arrow[from=1-2, to=1-3]
	\arrow["\wr", from=1-2, to=2-2]
	\arrow[from=1-3, to=1-4]
	\arrow[from=1-3, to=2-3]
	\arrow[from=1-4, to=1-5]
	\arrow["\wr", from=1-4, to=2-4]
	\arrow[from=2-1, to=2-2]
	\arrow[from=2-2, to=2-3]
	\arrow[from=2-3, to=2-4]
\end{tikzcd}\]
and the statement follows from the five lemma.

	\end{proof}
	\begin{corollary} \label{cor:vtiltingsoergelbimodules}
		Recall that the Kac--Moody datum $\Dd$ is cofree and of adjoint type. The functor $\V$ yields an monoidal equivalence between tilting sheaves in the universal monodromic Hecke category and $K$-theory Soergel bimodules 
		$$\Hecketilt\to \SBimK_{\widehat{\Dd}}$$
		mapping $T_s$ to $\LDR\otimes_{\LDR^s}\LDR.$
	\end{corollary}

\subsection{Soergel-theoretic description of the Hecke category}
	Combining the Struktursatz, see \Cref{thmStrukturSatzMonodromic} and \Cref{cor:vtiltingsoergelbimodules}, and the formality result, see \Cref{corGenerationBottSamelsonMonodromic}, we obtain the following `combinatorial' description of the universal monodromic Hecke category.
	\begin{theorem}\label{thm:soergelmonodromicheckecategory} Assume that the Kac--Moody datum $\Dd$ is cofree and of adjoint type. There is an equivalence of monoidal categories 
	\begin{equation*}
	\Heckemon_{\Dd} \rightarrow \Ch^b(\SBimK_{\widehat{\Dd}})
	\end{equation*}
	between the universal monodromic Hecke category and the category of bounded chain complexes of $K$-theoretic Soergel bimodules associated to $\widehat{\mathcal{D}}.$ 
	\end{theorem}

\section{Universal Koszul Duality}\label{sec:maintheorem}
Let $\Dd$ be a free Kac--Moody datum of simply-connected type. Recall that the Langlands dual Kac--Moody datum $\widehat{\Dd}$ is then cofree and of adjoint type, see \Cref{sec:kacmoodyrootdatum}.
We denote by $G\supset B\supset T$ the Kac--Moody group with Borel subgroup and maximal torus associated to $\Dd$ and $\widehat{G}\supset \widehat{B}\supset \widehat{T}$ be the Langlands dual groups associated to $\widehat{\Dd}$, see \Cref{sec:kacmoodygroups}. We denote by $\widehat{U}\subset \widehat{B}$ the unipotent radical.

In \Cref{sec:soergelktheoretichecke}, we studied the \emph{$K$-theoretic Hecke category} $\HeckeK_{\Dd}$ which consists of reduced $K$-motives on the stack $B\bs G/B.$
In \Cref{sec:soergelmonodromic}, we studied the \emph{universal monodromic Hecke category} $\Heckemon_{\widehat{\Dd}}$ which consists of $\widehat{T}\times \widehat{T}$-monodromic sheaves on the stack $\widehat{U}\bs \widehat{G}/\widehat{U}$. We note that in \Cref{sec:soergelmonodromic}, for notational ease, we worked with a group $G$ associated to a Kac--Moody datum $\Dd.$

By combining the Soergel-theoretic descriptions in terms of $K$-theory Soergel bimodules of both categories, see \Cref{thm:soergelktheoreticheckecategory} and \Cref{thm:soergelmonodromicheckecategory}, we obtain our main result.
\begin{theorem}[Universal Koszul Duality]\label{thm:main} There is a monoidal equivalence
    $$\HeckeK_\Dd\simeq \Heckemon_{\widehat\Dd}$$
which sends the pure $K$-motive $E_s$ to the tilting sheaf $T_s.$
\end{theorem}

\appendix
\section{Weight Structures}
We recall some basic properties of weight structures and weight complex functors.
\begin{definition}\cite[Definition 1.1.1]{bondarkoWeightStructuresVs2010}
	Let $\Cc$ be a triangulated category. A \emph{weight structure} $w$ on $\Cc$ is a pair $w=(\Cc^{w\leq 0},\Cc^{w\geq 0})$ of full subcategories of $\Cc,$ which are closed under direct summands, such that with $\Cc^{w\leq n}:=\Cc^{w\leq 0}[-n]$ and $\Cc^{w\geq n}:=\Cc^{w\geq 0}[-n]$ the following conditions are satisfied:
	\begin{enumerate}
		
		\item $\Cc^{w\leq 0}\subseteq \Cc^{w\leq 1}$ and $\Cc^{w\geq 1}\subseteq \Cc^{w\geq 0};$
		\item for all $X\in \Cc^{w\geq 0}$ and $Y\in\Cc^{w\leq -1}$, we have $\Hom{\Cc}{X}{Y}=0;$
		\item for any $X\in \Cc$ there is a distinguished triangle 
		\begin{center}\disttriangle{A}{X}{B}\end{center} with $A\in \Cc^{w\geq 1}$ and $B\in \Cc^{w\leq 0}.$
	\end{enumerate}
	The full subcategory $\Cc^{w=0}=\Cc^{w\leq 0}\cap\Cc^{w\geq 0}$ is called the \emph{heart of the weight structure}. The weight structure is called \emph{bounded} if $\Cc=\bigcup_n\Cc^{w\geq n}=\bigcup_n\Cc^{w\leq n}.$
\end{definition}
A weight structure on a stable $\infty$-category $\Cc$ is simply a weight structure on its homotopy category $\ho\Cc.$
Weight structures can be generated from their heart.
\begin{proposition}\label{prop:weightstructuresviaheart}
     Let $\Cc$ be a triangulated category and $\Tt\subset\Cc$ an additive idempotent-closed subcategory generating $\Cc$ as a triangulated category. Moreover, assume that $\Hom_\Cc(X,Y[n])=0$ for all $X,Y\in \Tt$ and $n>0.$ Then there is a unique weight structure on $\Cc$ such that $\Tt\subset \Cc^{w=0}$. Moreover, in this case $\Tt=\Cc^{w=0}.$
\end{proposition}
\begin{proof}
    This follows from \cite[Theorem 4.3.2]{bondarkoWeightStructuresVs2010}.
\end{proof}
There is the following dual construction to Beilinson's realization functor for $t$-structures, see \cite{beilinsonFaisceauxPervers1982}.
\begin{proposition}
    Let $\Cc$ be a stable $\infty$-category with a bounded weight structure, then there is a \emph{weight complex functor}
    $$\Cc\to \Ch^b(\ho\Cc^{w=0}).$$
    The weight complex functor commutes with functors that preserve the weight structure. Moreover, if $\Cc$ has a monoidal structure that preserves the weight heart, then the weight complex functor is monoidal.
\end{proposition}
\begin{proof}
    The compatibility with weight exact functors is shown in \cite{sosniloTheoremHeartNegative2017}. 
	In \cite{aokiWeightComplexFunctor2020} it is shown that the weight complex functor is \emph{symmetric} monoidal. A similar argument shows that it is also monoidal.
\end{proof}
We record the following special case in which the weight complex functor is an equivalence of categories.
\begin{proposition} \label{prop:weightcomplexequivalence}
    Let $\Cc$ be a stable $\infty$-category and $\Tt\subset\Cc$ an additive idempotent-closed subcategory generating $\Cc$ as a stable subcategory, such that $\Hom_\Cc(X,Y[n])=0$ for all $X,Y\in \Tt$ and $n\neq 0.$ Then the weight complex functor for the weight structure associated to $\Tt$ yields an equivalence of categories 
    $$\Cc\stackrel{\sim}{\to} \Ch^b(\ho\Tt).$$
\end{proposition}

\printbibliography

\end{document}